# THE PSEUDO-MARGINAL APPROACH FOR EFFICIENT MONTE CARLO COMPUTATIONS

BY CHRISTOPHE ANDRIEU[1] AND GARETH O. ROBERTS

*University of Bristol and University of Warwick*

We introduce a powerful and flexible MCMC algorithm for stochastic simulation. The method builds on a pseudo-marginal method originally introduced in [*Genetics* **164** (2003) 1139–1160], showing how algorithms which are approximations to an idealized marginal algorithm, can share the same marginal stationary distribution as the idealized method. Theoretical results are given describing the convergence properties of the proposed method, and simple numerical examples are given to illustrate the promising empirical characteristics of the technique. Interesting comparisons with a more obvious, but inexact, Monte Carlo approximation to the marginal algorithm, are also given.

**1. Introduction.** We are interested in the problem of simulation from a probability distribution $\pi(d\theta, dz)$ which, for now, we shall assume admits a density $\pi(\theta, z)$ with respect to some $\sigma$-finite measure (which we shall just write as $d\theta \times dz$). The variables $\theta$ and $z$ are elements of essentially arbitrary spaces, $\Theta$ and $\mathsf{Z}$, respectively. We partition the state space in this way because, either:

1. interest lies mainly in the marginal law $\pi(d\theta)$ of the variable $\theta \in \Theta$ [which we shall assume, for now, admits density $\pi(\theta)$ with respect to $d\theta$]; or
2. exploration of $\pi(\theta)$ by MCMC methods is more convenient by appropriate auxiliary simulation.

In a Bayesian framework, for example, $\theta$ could represent a parameter of interest and $z$ a set of missing data or latent variables; this includes, among others, hidden Markov models and their continuous generalizations, but also mixture models, and as we shall see in the application section,

Received March 2007; revised November 2007.
[1]Supported by an EPSRC Advance Research Fellowship.

model selection problems in general [8]. Often the variable $z$ is introduced for convenience, in particular in cases where the marginal density $\pi(\theta)$ is of sole interest. Indeed $\pi(\theta)$, or expectations with respect to it, might be analytically intractable or too complex to evaluate, whereas the introduction of $z$ might lead to an analytical expression, or ease the implementation of numerical methods.

A relatively generic way of numerically approximating such expectations consists of simulating an ergodic Markov chain $\{(\theta_i, z_i)\}$ which admits $\pi(\theta, z)$ as invariant probability density: such techniques are known under the acronym MCMC (Markov chain Monte Carlo). A typical sampling scheme will alternate sampling from the conditionals $\pi(\theta|z)$ and $\pi(z|\theta)$, or more generally ergodic Markov transition probabilities with these conditionals as invariant distributions. Although such so-called data augmentation schemes can very often ease programming and lead to elegant algorithms, it is well established that in numerous situations they can result in strongly positively correlated samples $\{(\theta_i, z_i)\}$ (see, e.g., [5, 9]), which is an undesirable property when efficiency is sought. On the other hand, if $\pi(\theta)$ was known analytically or cheap to compute, it would often be possible to generate "more efficient" samples $\{\theta_i\}$ from a Markov chain with a transition probability $P$, typically a Metropolis–Hastings (MH) transition with invariant density $\pi(\theta)$ and proposal density $q(\theta, \vartheta)$. This celebrated MCMC update consists, given that the Markov chain is currently at $\theta$, of proposing $\theta^* \sim q(\theta, \cdot)$ and set the next value of the chain $\vartheta = \theta^*$ with probability $\alpha(\theta, \theta^*)$ which depends on the values of the densities $\pi(\theta)$ and $\pi(\theta^*)$; otherwise we set $\vartheta = \theta$. More details are given below in the form of pseudo-code (note that, in our context, we shall term such an algorithm a *marginal* algorithm).

These general remarks have lead to the development of MCMC algorithms that try to combine the benefits of both approaches: possible statistical and computational efficiency of sampling directly from $\pi(\theta)$, and implementational ease of augmented, or auxiliary, schemes. A natural approach consists of approximating the intractable density values $\pi(\theta)$ and $\pi(\vartheta)$ required for the computation of the acceptance probability of the MH update with *importance sampling* estimates [8], that is for some integer $N \geq 1$ and some importance probability density $q_\theta(z)$ (satisfying the usual support assumption) one can consider the estimators

$$\tilde{\pi}^N(\theta) := \frac{1}{N} \sum_{k=1}^N \frac{\pi(\theta, z(k))}{q_\theta(z(k))} \qquad \text{with } z(k)|\theta \overset{\text{i.i.d.}}{\sim} q_\theta(\cdot)$$

and

$$(1.1) \qquad \tilde{\pi}^N(\vartheta) := \frac{1}{N} \sum_{k=1}^N \frac{\pi(\vartheta, \mathfrak{z}(k))}{q_\vartheta(\mathfrak{z}(k))} \qquad \text{with } \mathfrak{z}(k)|\vartheta \overset{\text{i.i.d.}}{\sim} q_\vartheta(\cdot),$$

TABLE 1
*Comparison of the marginal, MCWM and GIMH algorithms*

| Step | Marginal | MCWM | GIMH |
|---|---|---|---|
| 0. *Given:* | $\theta$ and $\pi(\theta)$ | $\theta$ and $\pi(\theta)$ | $\theta, Z$ and $\tilde{\pi}^N(\theta)$ |
| 1. *Sample:* | $\theta^* \sim q(\theta, \cdot)$ | $\theta^* \sim q(\theta, \cdot)$ | $\theta^* \sim q(\theta, \cdot)$ |
|  |  | $\begin{cases} Z \sim q_\theta^N(\cdot), \\ Z^* \sim q_{\theta^*}^N(\cdot) \end{cases}$ | $Z^* \sim q_{\theta^*}^N(\cdot)$ |
| 2. *Compute:* | $\pi(\theta^*)$ | $\begin{cases} \tilde{\pi}^N(\theta), \\ \tilde{\pi}^N(\theta^*) \end{cases}$ | $\tilde{\pi}^N(\theta^*)$ |
| 3. *Compute:* $r=$ | $\frac{\pi(\theta^*)q(\theta^*,\theta)}{\pi(\theta)q(\theta,\theta^*)}$ | $\frac{\tilde{\pi}^N(\theta^*)q(\theta^*,\theta)}{\tilde{\pi}^N(\theta)q(\theta,\theta^*)}$ | $\frac{\tilde{\pi}^N(\theta^*)q(\theta^*,\theta)}{\tilde{\pi}^N(\theta)q(\theta,\theta^*)}$ |
| 4. *With prob.* $1 \wedge r$: | $\vartheta = \theta^*$ | $\vartheta = \theta^*$ | $\begin{cases} \vartheta = \theta^*, \\ \mathfrak{z} = Z^* \end{cases}$ |
| *otherwise:* | $\vartheta = \theta$ | $\vartheta = \theta$ | $\begin{cases} \vartheta = \theta, \\ \mathfrak{z} = Z \end{cases}$ |

and simply plug these estimates in the expression for the marginal acceptance ratio (2.6). Note that here and hereafter we frequently omit the dependency on $Z := (z(1), z(2), \ldots, z(N))$ and $\mathfrak{Z} := (\mathfrak{z}(1), \mathfrak{z}(2), \ldots, \mathfrak{z}(N))$ for notational simplicity. We will denote $q_\theta^N(Z)$ and $q_\vartheta^N(\mathfrak{Z})$ the densities of $Z$ and $\mathfrak{Z}$.

There are, however, several possible implementations of this idea, and we now review two of them. Before embarking on a more formal presentation in Section 2, it can be helpful to give a comparative pseudo-code description of *MCWM* and *GIMH* (the acronyms are explained later in the text) and the *marginal* algorithm (see Table 1).

The first approach considered here, which corresponds to the middle column, is to attempt to approximate $P$, *independently at each iteration* using the importance ratio averages given in equation (1.1). More precisely both $Z$ and $Z^*$ are "refreshed" at each iteration independently of previously sampled auxiliary variables given $\theta$ and $\theta^*$. This algorithm is referred to as the Monte Carlo within Metropolis (*MCWM*) in [1] following the terminology of [7]. Due to the fact that the $Z$'s are independent at each iteration, one can easily see that $\{\theta_i\}$ is still a Markov chain with transition probability denoted $\tilde{P}_N^{\text{MCWM}}$ hereafter.

However MCWM and the marginal algorithm $P$ are not equivalent. In particular, $\pi(\theta)$ is typically not the invariant distribution density of $\tilde{P}_N^{\text{MCWM}}$ and therefore will not produce samples from $\pi(\theta)$ even in steady state. However, intuitively, provided that $\tilde{P}_N^{\text{MCWM}}$ is ergodic, the samples generated by this procedure will asymptotically be distributed according to an approximation of $\pi(\theta)$, which should be all the more precise that $N$ is large. Furthermore, we would like to know whether for sufficiently large $N$, the transition probability $\tilde{P}_N^{\text{MCWM}}$ does indeed inherit the convergence properties of $P$, as this was our initial motivation.

An interesting variation of *MCWM* could consist of using a single $Z$, sampled from some probability density $q_{\theta,\theta^*}^N(Z)$, to compute both $\tilde{\pi}^N(\theta)$ and $\tilde{\pi}^N(\theta^*)$. We do not pursue this here, but rather focus on the following.

In [1], Beaumont proposes a very interesting variation on the idea above, called grouped independence MH (*GIMH*), which corresponds to the rightmost column above. The MH transition probability of *GIMH*, as presented by Beaumont is similar in spirit to *MCWM*, but differs in that no fresh $Z$ is sampled at every iteration. Rather, *GIMH* can be interpreted as a form of *MCWM* where $Z$ is "recycled" from the previous iteration; as a result $Z$ is in general not distributed according to $q_\theta^N$, as is the case when the *MCWM* algorithm is used. Note that in addition $\{\theta_i\}$ is not a Markov chain anymore, but that $\{\theta_i, Z_i\}$ defines a Markov chain; we hereafter denote $\tilde{P}_N^{\text{GIMH}}$ the transition probability of this Markov chain.

The remarkable property noticed by Beaumont is that the acceptance ratio of *GIMH* can be rewritten as

$$\frac{\tilde{\pi}^N(\theta^*)q(\theta^*,\theta)}{\tilde{\pi}^N(\theta)q(\theta,\theta^*)}$$
$$= \frac{[1/N \sum_{k=1}^N \pi(\theta^*, z^*(k)) \prod_{l=1;l\neq k}^N q_{\theta^*}(z^*(l))]q(\theta^*,\theta)q_\theta^N(Z)}{[1/N \sum_{k=1}^N \pi(\theta, z(k)) \prod_{l=1;l\neq k}^N q_\theta(z(l))]q(\theta,\theta^*)q_{\theta^*}^N(Z^*)},$$

which suggests that $\tilde{P}_N^{\text{GIMH}}$ is, to complete Beaumont's argument, a MH algorithm with proposal density $q(\theta,\vartheta)q_\vartheta^N(\mathfrak{Z})$ and target density given between the brackets above, denoted $\tilde{\pi}^N(\theta, Z)$ hereafter. Hence, as soon as *GIMH* defines an irreducible and aperiodic Markov chain it will produce samples $\{\theta_i\}$ distributed in the limit as $i \to \infty$ according to the marginal $\pi(\theta)$. Given this interpretation, it is important to point out that updates of the type $\tilde{P}_N^{\text{GIMH}}$ are not the only available to sample from $\tilde{\pi}^N(\theta, Z)$, a point apparently missed in the literature. For example, it is possible to update "$Z$ given $\theta$," which is crucial to address some of the weaknesses of this approach, see Section 5. This distinction between target distribution and update also motivates the *pseudo-marginal* terminology adopted here.

The dual interpretation of *GIMH* as an approximation of a MH with target density $\pi(\theta)$ or an MH with target distribution density $\tilde{\pi}^N(\theta, Z)$ therefore opens the possibility for the design of algorithms that inherit the potential efficiency of $P$ while still being able to produce samples from $\pi(\theta)$, and not an approximation. However one can reiterate the questions asked earlier about the convergence properties of $\tilde{P}_N^{\text{MCWM}}$ and their relation to the ideal transition probability $P$.

Before giving some answers to these questions, we first show how the approach can be easily generalized in order to allow for more sophisticated transitions, leading to potential "local adaptation" schemes for example and also suggest new applications, such as model selection and applications to reversible jump MCMC algorithms [3]. Sections 3–5 are dedicated to the theoretical properties of generalizations of *GIMH*. Our main results are: Theorem 1, where we show that if the marginal chain is irreducible and aperiodic then generalizations of *GIMH* also converge; Theorem 6 shows that under very mild and intuitive conditions, mainly (A2) and (A3), generalizations of *GIMH* have finite horizon convergence properties very similar to those of the marginal algorithm, provided that $N$ is large enough; in Theorem 8, under more stringent assumption, we investigate the geometric and uniform convergence of generalizations of *GIMH*. Section 6 is dedicated to some theoretical properties of generalizations of *MCWM* which turn out to be much simpler to establish than for generalizations of *GIMH*. In particular we show in Theorem 9 that if the marginal algorithm is uniformly ergodic, then generalizations of *GIMH* can inherit this property with arbitrary precision. We conclude with Section 7 where we show how the ideas developed in this paper can be used in order to design efficient reversible jump MCMC algorithms to perform model selection using very simple mechanisms.

**2. Set up and notation.** Hereafter we will need the following notation. For some integer $N \geq 1$ let $Z := (z(1), z(2), \ldots, z(N)) \in \mathsf{Z}^N$ denote a generic vector of $\mathsf{Z}^N$ with coordinates $z(k)$, $k = 1, \ldots, N$. For any $Z \in \mathsf{Z}^N$ and $k = 1, \ldots, N$ we will denote $Z^{-k} := (z(1), \ldots, z(k-1), z(k+1), \ldots, z(N)) \in \mathsf{Z}^{N-1}$ with obvious conventions. For any $Z \in \mathsf{Z}^N$ and $\mathcal{Z} = (z(k_1), \ldots, z(k_l))$ (for $l \in \{1, \ldots, N-1\}$ and $k_1, \ldots, k_l \in \{1, \ldots, N\}$) a subvector of $Z$, we define $Z \setminus \{\mathcal{Z}\} := (z(1), \ldots, z(k_1 - 1), z(k_1 + 1), \ldots, z(k_2 - 1), z(k_2 + 1), \ldots)$ and $Z^l := (z(1), \ldots, z(l)) \in \mathsf{Z}^l$, with the notational convention $Z^0 = \varnothing$.

2.1. *The pseudo-marginal.* Let $(\Theta, \mathcal{B}(\Theta))$ and $(\mathsf{Z}, \mathcal{B}(\mathsf{Z}))$ be two measurable spaces. Let $\pi(d\theta, dz)$ be a probability distribution on the space $(\Theta \times \mathsf{Z}, \mathcal{B}(\Theta) \times \mathcal{B}(\mathsf{Z}))$, let $\pi(d\theta)$ be its marginal distribution and let us denote for any $\theta \in \Theta$, $\pi_\theta(dz)$ the associated conditional (on $\theta$) distribution. Let $\{Q_\theta^N(dZ), \theta \in \Theta \text{ and } N \in \mathbb{N}\}$ be a family of probability distributions, the "proposals," defined on $(\mathsf{Z}^N, \mathcal{B}(\mathsf{Z}^N))$, $\{(w_1^N, w_2^N, \ldots, w_N^N) \in [0,1]^N, \ N \in \mathbb{N}: \sum_{k=1}^N w_k^N = 1\}$ be a family of weights and let $\{\mathcal{Z}_k^N, k = 1, \ldots, N \text{ and } N \in \mathbb{N}\}$ be a family of subvectors of arbitrary sizes of vectors of the type $Z^{-k}$ as defined above. Before defining the pseudo-marginal and its associated joint model, we require the following assumption. Denoting for any $A \in \mathcal{B}(\mathsf{Z})$, $Q_\theta^N(z_k \in A | \mathcal{Z}_k^N)$ the conditional distribution or $z_k$ given $\mathcal{Z}_k^N$,

(A1) We assume that for our choice of $\{Q_\theta^N\}$, $\{w_k^N\}$ and $\{\mathcal{Z}_k^N\}$, for all $N \geq 1$, any $\theta \in \Theta$ and $k = 1, \ldots, N$, $\pi_\theta(\cdot) \ll Q_\theta^N(\cdot | \mathcal{Z}_k^N)$.

(A1) allows one to define for any $N \in \mathbb{N}$ the following linear combination of Radon–Nikodym derivatives for $(\theta, Z) \in \Theta \times \mathsf{Z}^N$ (the importance weights)

$$\gamma^N(\theta) := \sum_{k=1}^{N} w_k^N \frac{\pi_\theta(dz(k))}{Q_\theta^N(dz(k)|\mathcal{Z}_k^N)}, \tag{2.1}$$

the dependence on $Z$ being implicit. Whenever $\gamma^N(\theta) > 0$, we define $\lambda^N(\theta) := |\log \gamma^N(\theta)|$ and by convention we let $\lambda^N(\theta) := +\infty$ when $\gamma^N(\theta) = 0$. In turn we can define the following probability distribution on $(\Theta \times \mathsf{Z}, \mathcal{B}(\Theta) \times \mathcal{B}(\mathsf{Z}))$,

$$\tilde{\pi}^N(d\theta, dZ) := \pi(d\theta) Q_\theta^N(dZ) \gamma^N(\theta), \tag{2.2}$$

which, as we shall see, is a generalization of the underlying target distribution identified in equation (1.2). Hereafter for any $\theta \in \Theta$ we will denote

$$\tilde{\pi}_\theta^N(dZ) := Q_\theta^N(dZ) \gamma^N(\theta), \tag{2.3}$$

the *conditional probability distribution* of $Z$ given $\theta$ in equation (2.2) and we introduce for any $Z \in \mathsf{Z}$,

$$\tilde{\pi}^N(d\theta) := \pi(d\theta) \gamma^N(\theta)$$

the "*pseudo-marginal.*" Note that whenever $\pi(d\theta)$ has a density $\pi(\theta)$ and provided that $Z|\theta \sim Q_\theta^N$, then the associated probability density $\tilde{\pi}^N(\theta) := \pi(\theta) \gamma^N(\theta)$ is an unbiased importance sampling estimator of $\pi(\theta)$ a fundamental property which ensures that $\pi(d\theta)$ is the marginal of $\tilde{\pi}^N(d\theta, dZ)$. It can indeed be easily checked that for any $\theta \in \Theta$, $Q_\theta^N(\gamma^N(\theta)) = 1$. Note that in practice the variability of $\pi(\theta) \gamma^N(\theta)$ [and hence $Q_\theta^N(dZ)$] for any $\theta \in \Theta$ is expected to have an important influence on the performance of the algorithm an illustration is given in Theorem 8. We will frequently use the following identities:

$$\tilde{\pi}^N(d\theta, dZ) = \tilde{\pi}^N(d\theta) Q_\theta^N(dZ) = \pi(d\theta) \tilde{\pi}_\theta^N(dZ). \tag{2.4}$$

We conclude this section with various examples of choices of $\{w_i^N\}$ and $\{\mathcal{Z}_i^N\}$ introduced in the general framework presented earlier.

EXAMPLE 1 (Classical importance sampling). The case where $w_i^N = 1/N$ and $Q_\theta^N(dZ)$ is factorizable and exchangeable, that is,

$$Q_\theta^N(dZ) = \prod_{i=1}^{N} Q_\theta(dz(i)),$$

leads to Beaumont's *GIMH* algorithm.

EXAMPLE 2 (Sequential sampling). A choice of great practical interest, as illustrated later on in Section 7, consists of the case where $\mathcal{Z}_i^N = Z^{i-1}$, which allows for sequential sampling of $\{z(i)\}$, hence offering the possibility to adapt the sampling strategy in light of already sampled $z(i)$'s. This sequential framework encompasses the case where $\{z(i)\}$ is a realization from a Gibbs sampler with target distribution $\pi_\theta(dz)$ for some $\theta \in \Theta$. In such situations nondecreasing sequences $\{w_k^N\}$ might be preferable in order to discount the "burn-in" period.

EXAMPLE 3 (Gibbs sampler type). In some situations we might have good reasons to believe that the analytically intractable marginal distributions of $Q_\theta^N(dZ)$ for $\theta \in \Theta$ are good approximations of $\pi_\theta(dz)$. In this case one can suggest the application of the algorithm with $\mathcal{Z}_i^N = Z^{-i}$, which can be interpreted as a random scan Gibbs sampler to sample from $Q_\theta^N(dZ)$, and hence its marginals.

2.2. *Pseudo-marginal based algorithms.* We now introduce a formal description of the transition probabilities of the marginal algorithm and the two variants of the pseudo-marginal approach, which can be seen as generalization of *MCWM* and *GIMH*. The transition probability of the marginal algorithm, a standard MH algorithm, targets $\pi(d\theta)$ and uses $Q(\theta, d\vartheta)$ as proposal distribution is defined for any $\theta, \vartheta \in \Theta$ as

$$(2.5) \quad P(\theta, d\vartheta) := \alpha(\theta, \vartheta)\, Q(\theta, d\vartheta) + \delta_\theta(d\vartheta)\Big[1 - \int_\Theta \alpha(\theta, \vartheta) Q(\theta, d\vartheta)\Big],$$

where $\alpha(\theta, \vartheta) := 1 \wedge r(\theta, \vartheta)$ with (a)

$$(2.6) \quad 0 < r(\theta, \vartheta) := \frac{\pi(d\vartheta) Q(\vartheta, d\theta)}{\pi(d\theta) Q(\theta, d\vartheta)} < +\infty,$$

on a symmetric set $R \subset \Theta \times \Theta$ (see [11], Proposition 1 and Theorem 2), (b) $r(\theta, \vartheta) := 0$, $\pi(d\theta) Q(\theta, d\vartheta)$-almost everywhere on the complement $R^c$ of $R$ and (c) $r(\theta, \vartheta) := 1$ on measurable subsets of $R^c$ of $\pi(d\theta) Q(\theta, d\vartheta)$-zero probability.

The transition probability of the generalization of *MCWM* is not a standard MH algorithm. It consists of proposing $\vartheta \sim Q(\theta, \cdot)$, $Z \sim Q_\theta^N$ and $\mathfrak{Z} \sim Q_\vartheta^N$, compute an acceptance probability $\tilde{\alpha}^N(\theta, \vartheta)$ defined below, and accept or reject the proposal according to $\tilde{\alpha}^N(\theta, \vartheta)$. More formally the transition probability is defined for any $\theta, \vartheta \in \Theta$ as

$$(2.7)\quad \begin{aligned}\tilde{P}_N^{\text{noisy}}(\theta, d\vartheta) :=\ &Q_\theta^N \otimes Q_\vartheta^N(\tilde{\alpha}^N(\theta, \vartheta)) Q(\theta, d\vartheta) \\ &+ \delta_\theta(d\vartheta)\Big[1 - \int_\Theta Q_\theta^N \otimes Q_\vartheta^N(\tilde{\alpha}^N(\theta, \vartheta)) Q(\theta, d\vartheta)\Big],\end{aligned}$$

where $\otimes$ indicates the product for measures and $\tilde{\alpha}^N(\theta,\vartheta) := 1 \wedge \tilde{r}^N(\theta,\vartheta)$, with (a)

$$\tilde{r}^N(\theta,\vartheta) := \frac{\tilde{\pi}^N(d\vartheta)Q(\vartheta,d\theta)}{\tilde{\pi}^N(d\theta)Q(\theta,d\vartheta)}, \tag{2.8}$$

on an appropriate set $\tilde{R} \subset (\Theta \times \mathsf{Z}^N)^2$, (b) $\tilde{r}^N(\theta,\vartheta) := 0$, $\tilde{\pi}^N(d\theta)Q(\theta,d\vartheta)$-almost everywhere on the complement of $\tilde{R}$ and (c) 1 otherwise. Note that $\tilde{r}^N(\theta,\vartheta)$ can be computed even in situations where the normalizing constant of $\tilde{\pi}^N(d\theta, dZ)$ is unknown but that, on the other hand, the normalizing constant of $Q_\theta^N(dZ)$ might be required.

The transition probability of the *GIMH* variant of the pseudo-marginal approach is of the MH type and is defined on the extended space $\Theta \times \mathsf{Z}^N$. It targets $\tilde{\pi}^N(d\theta, dZ)$ and uses the proposal distribution $Q(\theta, d\vartheta)Q_\vartheta^N(d\mathfrak{Z})$

$$\begin{aligned}\tilde{P}_N^{\text{exact}}&(\theta, Z; d\vartheta, d\mathfrak{Z}) \\ &= \tilde{\alpha}^N(\theta,\vartheta) \, Q(\theta, d\vartheta)Q_\vartheta^N(d\mathfrak{Z}) \\ &\quad + \delta_{\theta,Z}(d\vartheta, d\mathfrak{Z})\Big[1 - \int_{\Theta \times \mathsf{Z}^N} \tilde{\alpha}^N(\theta,\vartheta)Q(\theta,d\vartheta)Q_\theta^N(d\mathfrak{Z})\Big],\end{aligned} \tag{2.9}$$

with $\tilde{\alpha}^N(\theta,\vartheta)$ as above equation (2.8). This expression for the acceptance probability of the exact pseudo-marginal algorithm relies on an identity, which we will frequently use later on, between the *marginal* acceptance ratio (2.6) and its exact *pseudo-marginal* counterpart,

$$\tilde{r}^N(\theta,\vartheta) := \frac{\tilde{\pi}^N(d\vartheta, d\mathfrak{Z})Q(\vartheta, d\theta)Q_\theta^N(dZ)}{\tilde{\pi}^N(d\theta, dZ)Q(\theta, d\vartheta)Q_\vartheta^N(d\mathfrak{Z})} = \frac{\gamma^N(\vartheta)}{\gamma^N(\theta)}r(\theta,\vartheta) \tag{2.10}$$

for any $(\theta, Z, \vartheta, \mathfrak{Z}) \in \tilde{R} := \{(\theta, Z, \vartheta, \mathfrak{Z}) : (\theta, \vartheta) \in R, \ Z \in \mathcal{Z}_\theta, \ \mathfrak{Z} \in \mathcal{Z}_\vartheta\}$ with $\mathcal{Z}_\theta := \{Z \in \mathsf{Z} : \gamma^N(\theta) > 0\}$. For any $N \in \mathbb{N}$ and $(\theta, Z) \in \Theta \times \mathsf{Z}^N$ we will denote $\alpha(\theta, Z)$ [resp. $\rho(\theta, Z)$]

$$\alpha(\theta, Z) := 1 - \rho(\theta, Z) := \int_{\Theta \times \mathsf{Z}^N} \tilde{\alpha}^N(\theta,\vartheta)Q(\theta, d\vartheta)Q_\vartheta^N(d\mathfrak{Z}), \tag{2.11}$$

the probability of leaving (resp. staying in) state $(\theta, Z)$. Note that we do not here make the dependence of this quantity on $N$ explicit for notational simplicity. Similarly we will denote $\alpha(\theta)$ [resp. $\rho(\theta)$] the probability of leaving (resp. staying in) state $\theta$ for transition $P$.

In the next two sections we study the Markov chain $\{\theta_i, Z_i\}$ started at $(\theta_0, Z_0) \in \Theta \times \mathsf{Z}^N$ and with transition probability $\tilde{P}_N^{\text{exact}}$ (which will be denoted $\tilde{P}_N$ for simplicity, when no ambiguity is possible) as given in equation (2.9) that is a MH with target distribution $\tilde{\pi}^N(d\theta, dZ)$ and proposal distribution $Q(\theta, d\vartheta)Q_\vartheta^N(d\mathfrak{Z})$ and acceptance ratio given by equation (2.10). In order to analyze the performance of the Markov chain generated by $\tilde{P}_N$ we will embed the exact marginal Markov chain with transition $P$ as in equation (2.5) defined on $(\Theta^{\mathbb{N}}, \mathcal{B}(\Theta^{\mathbb{N}}))$ into a Markov chain defined on $((\Theta \times \mathsf{Z}^N)^{\mathbb{N}}, (\mathcal{B}(\Theta) \times \mathcal{B}(\mathsf{Z}))^{\mathbb{N}})$ as follows. We define a Markov chain which is generated by a MH transition probability $\bar{P}_N$, with invariant distribution $\tilde{\pi}^N(d\theta, dZ)$ and proposal distribution $Q(\theta, d\vartheta) \tilde{\pi}_\vartheta^N(d\mathfrak{Z})$ [instead of $Q(\theta, d\vartheta) Q_\vartheta^N(d\mathfrak{Z})$ for $\tilde{P}_N$], leading to the transition probability,

$$\bar{P}_N(\theta, Z; d\vartheta, d\mathfrak{Z}) := \alpha(\theta, \vartheta) Q(\theta, d\vartheta) \tilde{\pi}_\vartheta^N(d\mathfrak{Z})$$
$$+ \delta_{\theta, Z}(d\vartheta, d\mathfrak{Z}) \left[1 - \int_\Theta \alpha(\theta, \vartheta) Q(\theta, d\vartheta)\right], \quad (2.12)$$

with $\alpha(\theta, \vartheta) := 1 \wedge r(\theta, \vartheta)$, where $r(\theta, \vartheta)$ is as in equation (2.6). Our analysis will rely upon a comparison of $\tilde{P}_N$ and $\bar{P}_N$, which live on a common space.

Finally, we will hereafter use the following standard notation for probabilities and Markov chain transition probabilities. For a space $(E, \mathcal{E})$ we define for $f : E \to \mathbb{R}$: for $\mu$ a measure on $(E, \mathcal{E})$ $\mu(f) := \int_E f(x) \mu(dx)$; $\|\mu\| := \frac{1}{2} \sup_{|f| \leq 1} |\mu(f)|$; for any $A \in \mathcal{E}$ $\mu(A) = \mu(\mathbb{I}\{x \in A\})$, where $\mathbb{I}\{\cdot \in A\} = \mathbb{I}\{A\}$ denotes the indicator function of set $A$; for a transition probability $\Pi : E \times \mathcal{E} \to [0, 1]$ $\Pi f(x) := \int_E \Pi(x, dy) f(y)$ and $\Pi^i f(x) := \Pi(x, \Pi^{i-1} f)$ for $i \geq 1$ with $\Pi^0 f(x) := f(x)$.

**3. A simple convergence result for exact algorithms.** The theory of $\psi$-irreducible Markov chains has proved to be a very powerful tool in order to analyze classical MCMC algorithms, and in particular the MH algorithm. More precisely, since MCMC deal with the situation where $\pi P = \pi$, then if in addition $P$ is $\psi$-irreducible and aperiodic, it can be shown that $\|P^k(\theta_0, \cdot) - \pi(\cdot)\| \to 0$ as $k \to \infty$ $\pi$-a.s. [4]. This motivates the following theorem.

THEOREM 1. *Assume* (A1) *and that $P$ defines a $\psi$-irreducible and aperiodic Markov chain such that $\pi P = \pi$. Then for any $N \in \mathbb{N}$ such that for any $(\theta, Z) \in \Theta \times \mathsf{Z}^N$, $\rho(\theta, Z) > 0$ [with $\rho(\theta, Z)$ as in equation (2.11)], $\tilde{P}_N$ is also $\psi$-irreducible and aperiodic, and hence $\tilde{\pi}^N$-a.s. [in $(\theta_0, Z_0) \in \Theta \times \mathsf{Z}^N$],*

$$\lim_{k \to +\infty} \|\tilde{P}_N^k(\theta_0, Z_0; \cdot) - \tilde{\pi}^N(\cdot)\| = 0.$$

PROOF. We here drop $N$ for simplicity. First notice that by construction if $P$ is $\psi$-irreducible and aperiodic, then so is $\bar{P}$ [defined in equation (2.12)] and consequently $\bar{P}$ defines an ergodic Markov chain with invariant distribution $\tilde{\pi}^N$. We will show that under the assumptions the accessible sets of $\bar{P}$ are included in those of $\tilde{P}$, which will allow us to conclude. More precisely we show by induction that for any $k \in \mathbb{N}$, $(\theta, Z) \in \Theta \times \mathsf{Z}^N$ and $A \times B \in$ $\mathcal{B}(\Theta) \times \mathcal{B}(\mathsf{Z}^N)$ such that $\bar{P}^k(\theta, Z; A \times B) > 0$, then $\tilde{P}^k(\theta, Z; A \times B) > 0$. For any $\theta \in \Theta$ recall that $\mathcal{Z}_\theta := \{Z : \gamma^N(\theta) > 0\}$ and for notational simplicity we will use the convention $1 \wedge \gamma^N(\vartheta)/\gamma^N(\theta) = 1$ whenever $\gamma^N(\theta) = 0$ below. First notice that from (2.10), for any $(\theta, Z) \in \Theta \times \mathsf{Z}^N$ and $A \times B \in \mathcal{B}(\Theta) \times \mathcal{B}(\mathsf{Z}^N)$,

$$
\begin{aligned}
\tilde{P}(\theta, Z; A \times B) &\geq \int_A Q_\vartheta^N \left(1 \wedge \frac{\gamma^N(\vartheta)}{\gamma^N(\theta)} \mathbb{I}(\mathfrak{Z} \in B \cap \mathcal{Z}_\vartheta)\right) \alpha(\theta, \vartheta) Q(\theta, d\vartheta) \\
&\quad + \mathbb{I}\{(\theta, Z) \in A \times B\} \rho(\theta, Z) \\
&\geq \int_A \tilde{\pi}_\vartheta^N \left(\frac{1}{\gamma^N(\vartheta)} 1 \wedge \frac{\gamma^N(\vartheta)}{\gamma^N(\theta)} \mathbb{I}(\mathfrak{Z} \in B \cap \mathcal{Z}_\vartheta)\right) \alpha(\theta, \vartheta) Q(\theta, d\vartheta) \\
&\quad + \mathbb{I}\{(\theta, Z) \in A \times B\} \rho(\theta, Z).
\end{aligned}
$$
(3.1)

Consequently, since for any $\theta \in \Theta$ and $B \in \mathcal{B}(\mathsf{Z}^N)$ we have $\tilde{\pi}_\theta^N(B) = \tilde{\pi}_\theta^N(B \cap \mathcal{Z}_\theta)$, we deduce that the implication is true for $k = 1$. Assume the induction assumption true up to some $k = n \geq 1$. Now for some $(\theta, Z) \in \Theta \times \mathsf{Z}^N$ let $A \times B \in \mathcal{B}(\Theta) \times \mathcal{B}(\mathsf{Z}^N)$ be such that $\bar{P}^{n+1}(\theta, Z; A \times B) > 0$ and assume that

$$\int_{\Theta \times \mathsf{Z}^N} \tilde{P}^n(\theta, Z; d\vartheta, d\mathfrak{Z}) \tilde{P}(\vartheta, \mathfrak{Z}; A \times B) = 0,$$

which implies that $\tilde{P}(\vartheta, \mathfrak{Z}; A \times B) = 0$, $\tilde{P}^n(\theta, Z; \cdot)$-a.s. and hence that $\bar{P}(\vartheta, \mathfrak{Z}; A \times B) = 0$, $\bar{P}^n(\theta, Z; \cdot)$-a.s. from the induction assumption for $k = 1$. From this and the induction assumption for $k = n$, we deduce that $\bar{P}(\vartheta, \mathfrak{Z}; A \times B) = 0$, $\bar{P}^n(\theta, Z; \cdot)$-a.s. (by contradiction), which contradicts the fact that $\bar{P}^{n+1}(\theta, Z; A \times B) > 0$. □

**4. Performance of the pseudo-marginal approach.** For the purpose of our analysis we introduce the following subsets of $\Theta$. For any $\epsilon > 0$, $N \in \mathbb{N}$ and denoting for any random variable $X$ with probability distribution $\mu$, $\mu(X > \epsilon) := \mu(\{X : X > \epsilon\})$,

(4.1) $\qquad \mathcal{T}(\epsilon, N) := \{\theta \in \Theta : Q_\theta^N(\lambda^N(\theta) > \epsilon) \leq \epsilon\},$

(4.2) $\qquad \mathcal{S}(\epsilon, N) := \{\theta \in \Theta : Q(\theta, \mathcal{T}(\epsilon, N)) \geq 1 - \epsilon\},$

(4.3) $\qquad \mathcal{R}(\epsilon, N) := \mathcal{S}(\epsilon, N) \cap \mathcal{T}(\epsilon, N),$

with $\lambda^N(\theta)$ as defined below equation (2.1). For a set $\mathcal{A} \subset \Theta$ we will denote $\bar{\mathcal{A}}$ its complement in $\Theta$.

The main result of this section is Theorem 6 and relies, in addition to (A1), on the following mild assumptions.

(A2) For any $\theta_0 \in \Theta$, $\lim_{k \to \infty} \|P^k(\theta_0, \cdot) - \pi(\cdot)\| = 0$.

(A3) For any $\theta \in \Theta$ and any $\epsilon > 0$,

$$\lim_{N \to \infty} Q_\theta^N(\lambda^N(\theta) > \epsilon) = 0.$$

Assumption (A3) is fundamental to our analysis, and implies the following two lemmata. First, it is a sufficient condition to control the total variation distance between $\tilde{\pi}_\theta^N$ and $Q_\theta^N$.

LEMMA 2. *Assume (A1) and (A3). Then for any $\theta \in \Theta$, $N \in \mathbb{N}$ and $\epsilon \in (0,1]$,*

$$\|\tilde{\pi}_\theta^N(\cdot) - Q_\theta^N(\cdot)\| \leq (3+e)\epsilon + 2Q_\theta^N(\lambda^N(\theta) > \epsilon)\mathbb{I}\{\theta \in \bar{\mathcal{T}}(\epsilon, N)\}.$$

PROOF. For any $\theta \in \Theta$, from equation (2.3)

$$\begin{aligned}(4.4) \quad \|\tilde{\pi}_\theta^N(\cdot) - Q_\theta^N(\cdot)\| &= Q_\theta^N(|\gamma^N(\theta) - 1|\mathbb{I}\{\lambda^N(\theta) > \epsilon\}) \\ &\quad + Q_\theta^N(|\gamma^N(\theta) - 1|\mathbb{I}\{\lambda^N(\theta) \leq \epsilon\}).\end{aligned}$$

On the one hand

$$(4.5) \qquad Q_\theta^N(|\gamma^N(\theta) - 1|\mathbb{I}\{\lambda^N(\theta) \leq \epsilon\}) \leq \exp(\epsilon) - 1 \leq e\epsilon,$$

and on the other hand

$$\begin{aligned}(4.6) \quad Q_\theta^N(|\gamma^N(\theta) - 1|\mathbb{I}\{\lambda^N(\theta) > \epsilon\}) &\leq Q_\theta^N(\gamma^N(\theta)\mathbb{I}\{\lambda^N(\theta) > \epsilon\}) \\ &\quad + Q_\theta^N(\lambda^N(\theta) > \epsilon).\end{aligned}$$

Notice that

$$(1-\epsilon)Q_\theta^N(\lambda^N(\theta) \leq \epsilon) \leq \exp(-\epsilon)Q_\theta^N(\lambda^N(\theta) \leq \epsilon) \leq Q_\theta^N(\gamma^N(\theta)\mathbb{I}\{\lambda^N(\theta) \leq \epsilon\}).$$

This, together with the fact that for any $\theta \in \Theta$, $Q_\theta^N(\gamma^N(\theta)) = \tilde{\pi}_\theta^N(1) = 1$, leads to

$$\begin{aligned}(4.7) \quad Q_\theta^N(\gamma^N(\theta)\mathbb{I}\{\lambda^N(\theta) > \epsilon\}) &= 1 - Q_\theta^N(\gamma^N(\theta)\mathbb{I}\{\lambda^N(\theta) \leq \epsilon\}) \\ &\leq 1 - (1-\epsilon)(1 - Q_\theta^N(\lambda^N(\theta) > \epsilon)) \\ &\leq \epsilon + Q_\theta^N(\lambda^N(\theta) > \epsilon).\end{aligned}$$

Now combining equations (4.4)–(4.7), we have for any $\theta \in \Theta$,

$$\begin{aligned}(4.8) \quad \|\tilde{\pi}_\theta^N(\cdot) - Q_\theta^N(\cdot)\| &\leq (1+e)\epsilon + 2Q_\theta^N(\lambda^N(\theta) > \epsilon) \\ &\leq (3+e)\epsilon + 2Q_\theta^N(\lambda^N(\theta) > \epsilon)\mathbb{I}\{\theta \in \bar{\mathcal{T}}(\epsilon, N)\}. \quad \square\end{aligned}$$

Assumption (A3) also implies the following important intermediate results.

LEMMA 3. *Let $\epsilon > 0$ and $\mathcal{S}(\epsilon, N), \mathcal{T}(\epsilon, N)$ and $\mathcal{R}(\epsilon, N)$ be as in (4.1)–(4.3). Assume* (A1) *and* (A3). *Then for any probability measure $\mu$ on $(\Theta, \mathcal{B}(\Theta))$,*

$$\lim_{N \to \infty} \mu(\mathcal{T}(\epsilon, N)) = 1, \tag{4.9}$$

$$\lim_{N \to \infty} \mu(\mathcal{S}(\epsilon, N)) = 1, \tag{4.10}$$

$$\lim_{N \to \infty} \mu(\mathcal{R}(\epsilon, N)) = 1. \tag{4.11}$$

PROOF. For any $\epsilon > 0$ and $\theta \in \Theta$, $\lim_{N \to \infty} \mathbb{I}\{\theta \in \mathcal{T}(\epsilon, N)\} = 1$ from (A3). Equation (4.9) follows from the dominated convergence theorem. To prove equation (4.10), note that equation (4.9) implies that for any $\theta \in \Theta$,

$$\lim_{N \to \infty} Q(\theta, \mathbb{I}\{\vartheta \in \mathcal{T}(\epsilon, N)\}) = 1.$$

Consequently for any $\theta \in \Theta$, $\lim_{N \to \infty} \mathbb{I}\{\theta \in \mathcal{S}(\epsilon, N)\} = 1$ and for any probability measure $\mu$, we have $\lim_{N \to \infty} \mu(\mathbb{I}\{\theta \in \mathcal{S}(\epsilon, N)\}) = \mu(1) = 1$. Equation (4.11) is immediate. $\square$

As we shall see, our results heavily rely on an estimate of the distance between $\bar{P}_N$ and $\tilde{P}_N$ under (A3).

LEMMA 4. *Assume* (A1) *and* (A3). *Let $\epsilon \in (0, 1]$ and $(\theta, Z) \in \mathcal{S}(\epsilon, N) \times \mathsf{Z}^N$ with $Z$ such that $\lambda^N(\theta) \leq \epsilon$ [$\mathcal{S}(\epsilon, N)$ being defined in (4.2)]. Then for $\bar{P}_N$ as defined in equation (2.12) and for any $\psi: \Theta \times \mathsf{Z}^N \to [-1, 1]$,*

$$|\bar{P}_N \psi(\theta, Z) - \tilde{P}_N \psi(\theta, Z)| \leq 24\epsilon.$$

PROOF. For simplicity we here drop $N$ in the notation for the transition probabilities. Let $(\theta, Z) \in \Theta \times \mathsf{Z}^N$ and $\psi: \Theta \times \mathsf{Z}^N \to [-1, 1]$. We have, by definition of $\bar{P}$,

$$\bar{P}\psi(\theta, Z) - \tilde{P}\psi(\theta, Z) = S_1 + S_2,$$

with

$$S_1 := \bar{P}\psi(\theta, Z) - \hat{P}\psi(\theta, Z),$$
$$S_2 := \hat{P}\psi(\theta, Z) - \tilde{P}\psi(\theta, Z),$$

where $\hat{P}$ is the MH transition probability with invariant distribution $\pi(d\theta)Q_\theta^N(dZ)$ and proposal distribution $Q(\theta, d\vartheta)Q_\vartheta^N(d\mathfrak{Z})$, that is

$$\hat{P}(\theta, Z; d\vartheta, d\mathfrak{Z}) := \alpha(\theta, \vartheta) Q(\theta, d\vartheta) Q_\vartheta^N(d\mathfrak{Z}) \tag{4.12}$$
$$+ \delta_{\theta, Z}(d\vartheta, d\mathfrak{Z}) \left[1 - \int_\Theta \alpha(\theta, \vartheta) Q(\theta, d\vartheta)\right].$$

From this and the definition of $\bar{P}$ in equation (2.12), the first term writes

(4.13)
$$|S_1| \leq |Q(\theta, (\tilde{\pi}_\vartheta^N - Q_\vartheta^N)\alpha(\theta,\vartheta)\psi)|$$
$$\leq 2Q(\theta, \|\tilde{\pi}_\vartheta^N(\cdot) - Q_\vartheta^N(\cdot)\|).$$

From this, Lemma 2 and since $\theta \in \mathcal{S}(\epsilon, N)$

(4.14)
$$|S_1| \leq 2[(3+e)\epsilon + 2Q(\theta, Q_\vartheta^N(\lambda^N(\vartheta) > \epsilon)\mathbb{I}\{\vartheta \in \bar{\mathcal{T}}(\epsilon, N)\})]$$
$$\leq 2(5+e)\epsilon.$$

The second term writes
$$S_2 = \int_{\Theta \times \mathsf{Z}^N} \psi(\vartheta, \mathfrak{z}) \Big[1 \wedge r(\theta, \vartheta) - 1 \wedge \frac{\gamma^N(\vartheta)}{\gamma^N(\theta)} r(\theta, \vartheta)\Big] Q(\theta, d\vartheta) Q_\vartheta^N(d\mathfrak{z})$$
$$- \psi(\theta, Z) \int_{\Theta \times \mathsf{Z}^N} \Big[1 \wedge r(\theta, \vartheta) - 1 \wedge \frac{\gamma^N(\vartheta)}{\gamma^N(\theta)} r(\theta, \vartheta)\Big] Q(\theta, d\vartheta) Q_\vartheta^N(d\mathfrak{z}),$$

and we therefore focus on the quantity
$$S_0 := \int_{\Theta \times \mathsf{Z}^N} \Big|1 \wedge r(\theta, \vartheta) - 1 \wedge \frac{\gamma^N(\vartheta)}{\gamma^N(\theta)} r(\theta, \vartheta)\Big| Q(\theta, d\vartheta) Q_\vartheta^N(d\mathfrak{z})$$
$$= \int_{\Theta \times \mathsf{Z}^N} \Big|1 \wedge r(\theta, \vartheta) - 1 \wedge \frac{\gamma^N(\vartheta)}{\gamma^N(\theta)} r(\theta, \vartheta)\Big| \mathbb{I}\{\lambda^N(\vartheta) > \epsilon\} Q(\theta, d\vartheta) Q_\vartheta^N(d\mathfrak{z})$$
$$+ \int_{\Theta \times \mathsf{Z}^N} \Big|1 \wedge r(\theta, \vartheta) - 1 \wedge \frac{\gamma^N(\vartheta)}{\gamma^N(\theta)} r(\theta, \vartheta)\Big| \mathbb{I}\{\lambda^N(\vartheta) \leq \epsilon\} Q(\theta, d\vartheta) Q_\vartheta^N(d\mathfrak{z}).$$

Noting that for any $(x,y) \in \mathbb{R}^2$,
$$|1 \wedge \exp(x) - 1 \wedge \exp(y)| = 1 \wedge |\exp(0 \wedge x) - \exp(0 \wedge y)| \leq 1 \wedge |x - y|,$$

we deduce that for $\theta \in \mathcal{S}(\epsilon, N)$,
$$S_0 \leq Q(\theta, Q_\vartheta^N(\mathbb{I}\{\lambda^N(\vartheta) > \epsilon\}))$$
$$+ Q(\theta, Q_\vartheta^N(1 \wedge |\log(\gamma^N(\theta)) - \log(\gamma^N(\vartheta))|\mathbb{I}\{\lambda^N(\vartheta) \leq \epsilon\})),$$

and consequently since $\gamma^N(\theta)$ depends on $\theta$ and $Z$ only,

(4.15)
$$|S_2| \leq 2(1 \wedge \lambda^N(\theta)) + 2Q(\theta, Q_\vartheta^N(\mathbb{I}\{\lambda^N(\vartheta) > \epsilon\}))$$
$$+ 2Q(\theta, Q_\vartheta^N(1 \wedge \lambda^N(\vartheta)\mathbb{I}\{\lambda^N(\vartheta) \leq \epsilon\})).$$

Consequently since $\theta \in \mathcal{S}(\epsilon, N)$,

(4.16)
$$|S_2| \leq 2(1 \wedge \lambda^N(\theta)) + 2Q(\theta, Q_\vartheta^N(\mathbb{I}\{\lambda^N(\vartheta) > \epsilon\}\mathbb{I}\{\vartheta \in \bar{\mathcal{T}}(\epsilon, N)\}))$$
$$+ 2Q(\theta, Q_\vartheta^N(1 \wedge \lambda^N(\vartheta)\mathbb{I}\{\lambda^N(\vartheta) \leq \epsilon\})\mathbb{I}\{\vartheta \in \bar{\mathcal{T}}(\epsilon, N)\})$$
$$+ 2Q(\theta, Q_\vartheta^N(\mathbb{I}\{\lambda^N(\vartheta) > \epsilon\}\mathbb{I}\{\vartheta \in \mathcal{T}(\epsilon, N)\}))$$
$$+ 2Q(\theta, Q_\vartheta^N(1 \wedge \lambda^N(\vartheta)\mathbb{I}\{\lambda^N(\vartheta) \leq \epsilon\})\mathbb{I}\{\vartheta \in \mathcal{T}(\epsilon, N)\})$$
$$\leq 2(\epsilon + \epsilon + \epsilon + \epsilon) = 8\epsilon.$$

One concludes by combining (4.14) and (4.16). □

We now combine Lemmata 3 and 4 to prove the following proposition.

PROPOSITION 5. *For any $\epsilon \in (0,1]$ and any probability measure $\mu$ on $(\Theta, \mathcal{B}(\Theta))$, there exists $N(\epsilon, \mu)$ such that for any $N \geq N(\epsilon, \mu)$ and any $\psi : \Theta \times \mathsf{Z}^N \to [-1, 1]$,*

$$|\mu(\tilde{\pi}_\theta^N((\bar{P}_N - \tilde{P}_N)\psi(\theta, Z)))| \leq \epsilon.$$

PROOF. Let $\varepsilon \in (0, 1]$, $\mu$ be a probability distribution on $(\Theta, \mathcal{B}(\Theta))$ and $\varphi : \Theta \to [-1, 1]$. We have, with $\mathcal{R}(\varepsilon, N)$ defined in equation (4.3),

$$|\mu(\varphi)| \leq |\mu(\varphi \mathbb{I}\{\theta \in \mathcal{R}(\varepsilon, N)\})| + \mu(\bar{\mathcal{R}}(\varepsilon, N)).$$

By Lemma 3 there exists $N_0(\varepsilon, \mu) \in \mathbb{N}$ (independent of $\varphi$) such that for $N \geq N_0(\varepsilon, \mu)$,

(4.17) $$|\mu(\bar{\mathcal{R}}(\varepsilon, N))| < \varepsilon.$$

For any $\psi : \Theta \times \mathsf{Z}^N \to [-1, 1]$, we have

(4.18) $$\mu(\mathbb{I}(\theta \in \mathcal{R}(\varepsilon, N))\tilde{\pi}_\theta^N((\tilde{P} - \bar{P})\psi)) = T_1 + T_2,$$

with

$$T_1 = \mu(\mathbb{I}(\theta \in \mathcal{R}(\varepsilon, N))\tilde{\pi}_\theta^N(\mathbb{I}(\lambda^N(\theta) \leq \varepsilon)(\tilde{P} - \bar{P})\psi)),$$
$$T_2 = \mu(\mathbb{I}(\theta \in \mathcal{R}(\varepsilon, N))\tilde{\pi}_\theta^N(\mathbb{I}(\lambda^N(\theta) > \varepsilon)(\tilde{P} - \bar{P})\psi)).$$

We apply Lemma 4 to $T_1$ and conclude that

(4.19) $$|T_1| \leq 24\varepsilon.$$

We now turn to $T_2$. First recall from equation (2.3) that for any $\theta \in \Theta$,

$$\tilde{\pi}_\theta^N(\lambda^N(\theta) \leq \varepsilon) = Q_\theta^N(\gamma^N(\theta)\mathbb{I}\{\lambda^N(\theta) \leq \varepsilon\}),$$

and, since $Q_\theta^N(\lambda^N(\theta) > \varepsilon) \leq \varepsilon$ for $\theta \in \mathcal{T}(\varepsilon, N)$, we conclude that for any $\theta \in \mathcal{T}(\varepsilon, N)$

$$(1 - \varepsilon)^2 \leq (1 - \varepsilon)Q_\theta^N(\lambda^N(\theta) \leq \varepsilon)$$
$$\leq Q_\theta^N(\exp(-\varepsilon)\mathbb{I}\{\lambda^N(\theta) \leq \varepsilon\}) \leq \tilde{\pi}_\theta^N(\lambda^N(\theta) \leq \varepsilon).$$

Hence, from the definition of $\mathcal{R}(\varepsilon, N)$,

$$|T_2| = |\mu(\mathbb{I}\{\theta \in \mathcal{R}(\varepsilon, N)\}\tilde{\pi}_\theta^N(\mathbb{I}\{\lambda^N(\theta) > \varepsilon\}(\tilde{P} - \bar{P})\psi))|$$
(4.20) $$\leq 2|\mu(\mathbb{I}\{\theta \in \mathcal{R}(\varepsilon, N)\})|(1 - (1 - \varepsilon)^2)$$
$$\leq 4\varepsilon.$$

Now we choose $\varepsilon = \epsilon/30$ and conclude with $N(\epsilon, \mu) = N_0(\varepsilon, \mu)$ and by combining (4.17) (for $2\varphi$), (4.19) and (4.20). □

Our main result is as follows and provides us with a bound $\ell$ on the loss of efficiency of the approximating chain compared to the ideal chain, which can be made arbitrarily small for large $N$'s. Define for any $\theta \in \Theta$ and any $\varepsilon \in (0, 1]$, $k(\varepsilon, \theta) := \inf\{k : \|P^k(\theta, \cdot) - \pi(\cdot)\| \leq \varepsilon\}$ and recall that for any $\theta \in \Theta$, $\rho(\theta) := 1 - \int_\Theta \alpha(\theta, \vartheta) Q(\theta, d\vartheta)$ is the probability of not leaving $\theta$ for $P$.

THEOREM 6. *Assume* (A1), (A2) *and* (A3). *Let* $\epsilon, \ell > 0$ *and* $\theta_0 \in \Theta$. *Then there exists* $N(\epsilon, \ell, \theta_0) \in \mathbb{N}$ *such that for any* $N \geq N(\epsilon, \ell, \theta_0)$ *and* $Z_0 \in \Theta \times \mathsf{Z}^N$ *such that* $\lambda^N(\theta_0) < \ell\epsilon/(24 k(\epsilon, \theta_0))$ *we have for any* $k \geq k(\epsilon, \theta_0)$,

$$\|\tilde{P}_N^k(\theta_0, Z_0; \cdot) - \tilde{\pi}^N(\cdot)\| \leq (1 + \ell)\epsilon + \rho^k(\theta_0).$$

COROLLARY 7. *Under the assumptions of Theorem 6, for any* $\epsilon, \ell > 0$ *and* $\theta_0 \in \Theta$, *there exists* $N(\epsilon, \ell, \theta_0) \in \mathbb{N}$ *such that for any* $N \geq N(\epsilon, \ell, \theta_0)$ *and* $Z_0 \in \Theta \times \mathsf{Z}^N$ *such that* $\lambda^N(\theta_0) < \ell\epsilon/(24 k(\epsilon, \theta_0))$ *we have for any* $k \geq k(\epsilon, \theta_0)$ *and any* $\varphi : \Theta \to [-1, 1]$,

$$\tfrac{1}{2}|\tilde{P}_N^k(\theta_0, Z_0; \varphi) - \pi(\varphi)| \leq (1 + \ell)\epsilon + \rho^k(\theta_0).$$

PROOF. Dropping $N$ for notational simplicity, we have that for any $k \geq 1$, $(\theta_0, Z_0) \in \Theta \times \mathsf{Z}^N$ and any $\psi : \Theta \times \mathsf{Z}^N \to [-1, 1]$,

$$\tilde{P}^k \psi(\theta_0, Z_0) - \tilde{\pi}(\psi) = S_0(k) + S_1(k) + S_2(k), \tag{4.21}$$

with $[\tilde{\pi}_\theta^N(\psi) := \tilde{\pi}_\theta^N(\psi(\theta, \cdot))$ hereafter for notational simplicity]

$$S_0(k) = \bar{P}^k \psi(\theta_0, Z_0) - P^k(\tilde{\pi}_\theta^N(\psi))(\theta_0),$$
$$S_1(k) = P^k(\tilde{\pi}_\theta^N(\psi))(\theta_0) - \pi(\tilde{\pi}_\theta^N(\psi)),$$
$$S_2(k) = \tilde{P}^k \psi(\theta_0, Z_0) - \bar{P}^k \psi(\theta_0, Z_0),$$

where the magnitude of $S_1(k)$ can be controlled thanks to the properties of the transition probability $P$, $S_0(k)$ and $S_2(k)$ correspond to the bias introduced by the approximation to the "ideal" chain. As we shall see, for a fixed $k$ this bias can be made arbitrarily small for $N$ sufficiently large. Let $\epsilon > 0$ and $(\theta_0, Z_0) \in \Theta \times \mathsf{Z}^N$ such that $\lambda^N(\theta_0) < \epsilon$. By a coupling argument or induction $|S_0(k)| \leq 2\, \rho(\theta_0)^k$. Since $\tilde{\pi}_\theta^N(\psi) : \Theta \to [-1, 1]$, by (A2) we have

$$k(\epsilon, \theta_0) < +\infty \quad \text{and} \quad |S_1(k(\epsilon, \theta_0))| \leq 2\epsilon. \tag{4.22}$$

From now on we set $k_0 := k(\epsilon, \theta_0)$ and use the following telescoping sum decomposition:

$$S_2 := S_2(k_0) = \sum_{l=0}^{k_0 - 1} \bar{P}^l \tilde{P}^{k_0 - l} \psi(\theta_0, Z_0) - \bar{P}^{l+1} \tilde{P}^{k_0 - (l+1)} \psi(\theta_0, Z_0)$$

$$= \sum_{l=0}^{k_0-1} \bar{P}^l (\tilde{P} - \bar{P}) \tilde{P}^{k_0-(l+1)} \psi(\theta_0, Z_0).$$

Let $\varepsilon \in (0,1]$. Noticing for any $l > 1$ we have for any $\bar{\psi}: \Theta \times \mathsf{Z}^N \to [-1,1]$

$$\bar{P}^l \bar{\psi}(\theta_0, Z_0) = \rho(\theta_0)^l \bar{\psi}(\theta_0, Z_0)$$
$$+ \sum_{j=1}^{l} \bar{P}^{j-1} \{Q(\theta_{j-1}, \alpha(\theta_{j-1}, \theta_j) \rho(\theta_j)^{l-j} \tilde{\pi}_{\theta_j}^N (\bar{\psi}(\theta_j, \cdot)))\}(\theta_0, Z_0)$$

we apply Lemma 4 $k_0$ times and Proposition 5 $(k_0 - 1)$ times (the result trivially applies to any finite measure) to show that there exists $N(\varepsilon, \theta_0)$ such that for any $N \geq N(\varepsilon, \theta_0)$ and some $C < \ell/(24k_0)$

(4.23) $$|S_2| \leq 24 C k_0 \epsilon + (k_0 - 1) \varepsilon.$$

We conclude by taking $\varepsilon = 2\epsilon(\ell - 24 k_0 C)/(k_0 - 1)$ in equation (4.23) and combining with equation (4.22) in equation (4.21). □

**5. Uniform and geometric ergodicity of exact algorithms.** In this section we illustrate the critical importance of the choice of a good importance sampling distribution $Q_\theta^N$ to ensure that $\tilde{P}_N$ is uniformly ergodic. More precisely we show that, for a given $N \in \mathbb{N}$, if the importance weights $\gamma^N(\theta)$ involved in the definition of the pseudo-marginal $\tilde{\pi}^N(d\theta)$ are unbounded for "too many" $\theta$'s, then $\tilde{P}_N$ cannot be geometrically ergodic. As we shall see "too many" will be quantified in terms of the measure of the set $\mathcal{U}^N := \{\theta: \forall M > 0, \; Q_\theta^N(\gamma^N(\theta) > M) > 0\}$ under $\pi$. In addition, for a fixed $N \in \mathbb{N}$, using the fact that it is most often possible to prove the uniform ergodicity of the MH update $P$ defined in (2.5) by establishing a minorization condition for the sub-stochastic kernel $K(\theta, d\vartheta) := \alpha(\theta, \vartheta) Q(\theta, d\vartheta)$ [that is there exists $n_0 \geq 1$, a constant $\epsilon > 0$ and a probability measure $\nu$ on $(\Theta, \mathcal{B}(\Theta))$ such that for any $(\theta, A) \in \Theta \times \mathcal{B}(\Theta)$, $K^{n_0}(\theta, A) \geq \epsilon \nu(A)$], we show that this property is systematically inherited by $\tilde{P}_N$ whenever $\gamma_*^N := \sup_{\theta, Z \in \Theta \times \mathsf{Z}^N} \gamma^N(\theta) < +\infty$.

THEOREM 8. *Assume* (A1) *and let* $N \in \mathbb{N}$. *Then:*

1. *if* $\pi(\mathcal{U}^N) > 0$, *then* $\tilde{P}_N$ *cannot be geometrically ergodic;*
2. *if we assume that there exist* $n_0 \geq 1$, *a constant* $\epsilon > 0$ *and a measure* $\nu$ *on* $(\Theta, \mathcal{B}(\Theta))$ *such that for any* $\theta, A \in \Theta \times \mathcal{B}(\Theta)$, $K^{n_0}(\theta, A) \geq \epsilon \nu(A)$ *(which implies that* $P$ *is uniformly ergodic) then if in addition* $\gamma_*^N < +\infty$ *then* $\tilde{P}_N$ *is uniformly ergodic.*

REMARK 1. Concerning the second point of the theorem, it should be noted that it is not possible in general to achieve the rate of convergence of the marginal chain $P$, even when $\{\gamma_*^N\}$ is bounded. Indeed, consider the independent MH algorithm, in the discrete case for simplicity and densities with respect to the counting measure. It is possible to characterize exactly the second-largest eigenvalue of the transition probability. For $P$ it takes the form $1 - (\sup_{\theta \in \Theta} \frac{\pi(\theta)}{q(\theta)})^{-1}$, while for Beaumont's form of the pseudo-marginal algorithm it will take the form $1 - (\sup_{(\theta,z) \in \Theta \times \mathsf{Z}} \frac{\pi(\theta,z)}{q(\theta,z)})^{-1}$. In the particular case where $\pi(\theta, z) = \pi(\theta)\pi(z)$ and $q(\theta, z) = q(\theta)q(z)$ this latter expression becomes $1 - (\sup_{\theta \in \Theta} \frac{\pi(\theta)}{q(\theta)})^{-1}(\sup_{z \in \mathsf{Z}} \frac{\pi(z)}{q(z)})^{-1}$ which in general will be larger than $1 - (\sup_{\theta \in \Theta} \frac{\pi(\theta)}{q(\theta)})^{-1}$. As we shall see this is not the case for the MCWM algorithm under appropriate assumptions (see Theorem 9).

PROOF OF THEOREM 8. We drop $N$ in $\tilde{P}_N$ and prove the first statement. We want to show that under the stated assumptions, for any $\epsilon > 0$

$$\tilde{\pi}^N(\mathbb{I}\{\alpha(\theta, Z) \leq \epsilon\}) = \pi\{Q_\theta^N(\gamma^N(\theta)\mathbb{I}\{\alpha(\theta, Z) \leq \epsilon\})\} > 0,$$

where $\alpha(\theta, Z)$ is defined in equation (2.11). From [10], Theorem 5.1, this indeed implies that $\tilde{P}$ cannot be geometric. For any $(\theta, Z) \in \Theta \times \mathsf{Z}^N$ with $\gamma^N(\theta) > 0$, by Fubini's theorem, Jensen's inequality and since for any $\vartheta \in \Theta$, $Q_\vartheta^N(\gamma^N(\vartheta)) = 1$, we have

$$\alpha^N(\theta, Z) = \int_{\Theta \times \mathsf{Z}^N} \left(1 \wedge \frac{\gamma^N(\vartheta)}{\gamma^N(\theta)} r(\theta, \vartheta)\right) Q(\theta, d\vartheta) Q_\vartheta^N(d\mathfrak{Z})$$

$$\leq \int_\Theta \left(1 \wedge \frac{r(\theta, \vartheta)}{\gamma^N(\theta)}\right) Q(\theta, d\vartheta).$$

For any $(\theta, \vartheta) \in \mathcal{U}^N \times \Theta$, since $r(\theta, \vartheta) < +\infty$ [see equation (2.6)] and

$$\{Z : \gamma^N(\theta) > M\} \neq \varnothing$$

for any $M > 0$, we have

$$\lim_{M \to \infty} \sup_{\{Z : \gamma^N(\theta) > M\}} 1 \wedge \frac{r(\theta, \vartheta)}{\gamma^N(\theta)} = 0.$$

Consequently, by the dominated convergence theorem,

$$\lim_{M \to \infty} \sup_{\{Z : \gamma^N(\theta) > M\}} \alpha(\theta, Z) = 0,$$

hence for any $\epsilon > 0$ there exists $M < +\infty$ such that

$$Q_\theta^N(\{Z : \gamma^N(\theta) > M, \alpha(\theta, Z) \leq \epsilon\}) > 0$$

and hence

$$Q_\theta^N(\gamma^N(\theta)\mathbb{I}\{\alpha(\theta, Z) \leq \epsilon\}) \geq M Q_\theta^N(\mathbb{I}\{\gamma^N(\theta) > M, \ \alpha(\theta, Z) \leq \epsilon\}) > 0.$$

We deduce that

$$\pi\{Q_\theta^N(\gamma^N(\theta)\mathbb{I}\{\alpha(\theta,Z)\leq\epsilon\})\}>0.$$

We now turn to the proof of the second claim. We first show by induction that for any $k\geq 1$ and $A\times B\in\mathcal{B}(\Theta)\times\mathcal{B}(\mathsf{Z}^N)$, $\tilde{P}^k(\theta,Z;A\times B)\geq\gamma_*^{-k}\int_A K^k(\theta,d\vartheta)\tilde{\pi}_\vartheta^N(B)$. For $k=1$,

$$\begin{aligned}
\tilde{P}(\theta,Z;A\times B) &\geq \int_{A\times B}\tilde{\alpha}_N(\theta,\vartheta)Q(\theta,d\vartheta)Q_\vartheta(d\mathfrak{Z}) \\
&\geq \int_{A\times B}\left(1\wedge\frac{\gamma^N(\vartheta)}{\gamma^N(\theta)}\right)\alpha(\theta,\vartheta)Q(\theta,d\vartheta)Q_\vartheta^N(d\mathfrak{Z}) \\
&\geq \gamma_*^{-1}\int_{A\times B}Q_\vartheta^N(\gamma^N(\vartheta))\alpha(\theta,\vartheta)Q(\theta,d\vartheta) \\
&= \gamma_*^{-1}\int_A K(\theta;d\vartheta)\tilde{\pi}_\vartheta^N(B).
\end{aligned}$$

(5.1)

Assume the inequality is true for some $k\geq 1$. Then from the induction assumption and equation (5.1),

$$\begin{aligned}
\tilde{P}^{k+1}(\theta,Z;A\times B) &= \int_{\Theta\times\mathsf{Z}^N}\tilde{P}^k(\theta,Z;d\vartheta,d\mathfrak{Z})\tilde{P}(\vartheta,\mathfrak{Z};A\times B) \\
&\geq \int_\Theta \gamma_*^{-k}K^k(\theta,d\vartheta)\tilde{\pi}_\vartheta^N(\mathsf{Z}^N)\gamma_*^{-1}\int_A K(\theta;d\vartheta)\tilde{\pi}_\vartheta^N(B) \\
&\geq \gamma_*^{-(k+1)}\int_A K^{k+1}(\theta,d\vartheta)\tilde{\pi}_\vartheta^N(B).
\end{aligned}$$

Hence the result. From this result for $k=n_0$ and the minorization assumption on $K$ we deduce that for any $(\theta,Z,A\times B)\in\Theta\times\mathsf{Z}^N\times(\mathcal{B}(\Theta)\times\mathcal{B}(\mathsf{Z}^N))$

$$\begin{aligned}
\tilde{P}^{n_0}(\theta,Z;A\times B) &\geq \gamma_*^{-n_0}\int_A K^{n_0}(\theta,d\vartheta)\tilde{\pi}_\vartheta^N(B) \\
&\geq \epsilon\gamma_*^{-n_0}\nu(\mathbb{I}_A\{\vartheta\}\tilde{\pi}_\vartheta^N(B)),
\end{aligned}$$

and hence the second claim. □

**6. Epilogue: MCWM.** As pointed out earlier, the analysis of generalizations of *MCWM*, defined in equation (2.7), is simpler than that of *GIMH* generalizations and relies on more classical arguments. This is due mainly to the fact that in this case $Z\sim Q_\theta^N$ and $\mathfrak{Z}\sim Q_\vartheta^N$. However the existence of an invariant distribution for $\tilde{P}_N^{\text{noisy}}$ is not obvious in general (it is not a MH update). This is in contrast with $\tilde{P}_N^{\text{exact}}$, for which the invariant distribution as well as its marginal distribution are known to be $\tilde{\pi}^N(d\theta,dZ)$ and $\pi(d\theta)$, respectively. We give here a result which characterizes the invariant $\check{\pi}^N(d\theta)$ when it exists, and the rate of convergence of $\tilde{P}_N^{\text{noisy}}$ (denoted $\tilde{P}_N$ in this section for simplicity), when $P$ is uniformly ergodic, and a simple uniform weak law of large numbers holds for $\lambda^N(\theta)$.

(A4) There exist $C \in (0, +\infty)$ and $\rho \in (0,1)$ such that for any $\theta_0 \in \Theta$ and $k \in \mathbb{N}$

$$\|P^k(\theta_0, \cdot) - \pi(\cdot)\| \leq C\rho^k.$$

(A5) We assume that for any $\epsilon > 0$,

$$\lim_{N \to \infty} \sup_{\theta \in \Theta} Q_\theta^N(\lambda^N(\theta) > \epsilon) = 0.$$

THEOREM 9. *Assume* (A1), (A4), (A5) *and that for any $N \geq 1$ there exists a probability distribution $\check{\pi}^N$ on $(\Theta, \mathcal{B}(\Theta))$ such that $\check{\pi}^N \tilde{P}_N = \check{\pi}^N$, with $\tilde{P}_N = \tilde{P}_N^{\text{noisy}}$ defined in equation (2.7). Then for any $\epsilon \in (0, \rho^{-1} - 1)$ there exists $N(\epsilon, \rho) \in \mathbb{N}$, $\tilde{\rho} \in (\rho, \rho(1+\epsilon)] \subset (\rho, 1)$ and $\tilde{C} \in (0, +\infty)$ such that for all $N \geq N(\epsilon, \rho)$, $\theta_0 \in \Theta$ and $k \geq 1$,*

(6.1) $$\|\tilde{P}_N^k(\theta_0, \cdot) - \check{\pi}^N(\cdot)\| \leq \tilde{C}\tilde{\rho}^k,$$

(6.2) $$\|\pi - \check{\pi}^N\| \leq C \frac{\epsilon}{1 - \rho}.$$

PROOF. In some instances we here drop $N$ for notational simplicity. First notice that for any $\theta \in \Theta$ and $n \in \mathbb{N}$,

$$\|\tilde{P}^n(\theta, \cdot) - P^n(\theta, \cdot)\| \leq \sum_{i=0}^{n-1} \|\tilde{P}^{n-i-1}(\theta, (P - \tilde{P})(P^i - \pi)(\cdot))\|$$

(6.3)
$$\leq C \sup_{\theta \in \Theta} \|P(\theta, \cdot) - \tilde{P}(\theta, \cdot)\| \sum_{i=0}^{n-1} \rho^i$$

$$\leq \frac{C}{1 - \rho} \sup_{\theta \in \Theta} \|P(\theta, \cdot) - \tilde{P}(\theta, \cdot)\|.$$

We bound the last term on the right-hand side. We have

$$P(\theta, d\vartheta) - \tilde{P}(\theta, d\vartheta) = Q_\theta^N \otimes Q_\vartheta^N(\alpha(\theta, \vartheta) - \tilde{\alpha}(\theta, Z; \vartheta, \mathfrak{Z}))Q(\theta, d\vartheta)$$

$$+ \delta_\theta(d\vartheta) \int_\Theta Q_\theta^N \otimes Q_\vartheta^N(\tilde{\alpha}(\theta, Z; \vartheta, \mathfrak{Z}) - \alpha(\theta, \vartheta))$$

$$\times Q(\theta, d\vartheta).$$

Let $\varepsilon \in (0, 1]$ and notice that since

$$\{Z, \mathfrak{Z} \in \mathsf{Z}^N\} \subset \{Z, \mathfrak{Z} : \lambda^N(\theta) > \varepsilon\} \cup \{Z, \mathfrak{Z} : \lambda^N(\vartheta) > \varepsilon\}$$
$$\cup \{Z, \mathfrak{Z} : \lambda^N(\theta) \leq \varepsilon, \lambda^N(\vartheta) \leq \varepsilon\},$$

we have

$$|Q_\theta^N \otimes Q_\vartheta^N(\alpha(\theta,\vartheta) - \tilde{\alpha}(\theta, Z; \vartheta, \mathfrak{Z}))| \leq Q_\theta^N(\lambda^N(\theta) > \varepsilon) + Q_\vartheta^N(\lambda^N(\vartheta) > \varepsilon)$$
$$+ Q_\theta^N(1 \wedge \lambda^N(\theta)\mathbb{I}(\lambda^N(\theta) \leq \varepsilon))$$
$$+ Q_\vartheta^N(1 \wedge \lambda^N(\vartheta)\mathbb{I}(\lambda^N(\vartheta) \leq \varepsilon)).$$

Following the proof of Lemma 4, and from (A5), we conclude that there exists $N(\varepsilon)$ such that for $N \geq N(\varepsilon)$ and any $\theta, \vartheta \in \Theta$,

$$|Q_\theta^N \otimes Q_\vartheta^N(\alpha(\theta,\vartheta) - \tilde{\alpha}(\theta, Z; \vartheta, \mathfrak{Z}))| \leq 4\varepsilon.$$

Consequently for any $\varepsilon \in (0,1]$ there exists $N(\varepsilon) \in \mathbb{N}$ such that for any $N \geq N(\varepsilon)$ and $\theta \in \Theta$,

$$\|P(\theta, \cdot) - \tilde{P}(\theta, \cdot)\| \leq 4\varepsilon.$$

As a result and from (A4) for any $(\theta, \vartheta) \in \Theta$ and $n \in \mathbb{N}$,

$$\|\tilde{P}^n(\theta, \cdot) - \tilde{P}^n(\vartheta, \cdot)\| \leq \|\tilde{P}^n(\theta, \cdot) - P^n(\theta, \cdot)\| + \|P^n(\vartheta, \cdot) - \tilde{P}^n(\vartheta, \cdot)\|$$
$$+ \|P^n(\vartheta, \cdot) - P^n(\theta, \cdot)\|$$
$$\leq C\rho^n + \frac{8C\varepsilon}{1-\rho}.$$

Define

$$\tilde{\rho} := \rho \sqrt[n]{C\left(1 + \frac{8\varepsilon\rho^{-n}}{1-\rho}\right)} \leq \rho \sqrt[n]{C}\left(1 + \frac{1}{n}\frac{8\varepsilon\rho^{-n}}{1-\rho}\right).$$

Choose $\epsilon \in (0, \rho^{-1} - 1)$ and let $n \in \mathbb{N}$ be such that $\sqrt[n]{C} \leq \sqrt{1+\epsilon}$ and $\varepsilon$ (depending on $n$ and $\rho$) be such that

$$1 + \frac{1}{n}\frac{8\varepsilon\rho^{-n}}{(1-\rho)} \leq \sqrt{1+\epsilon}.$$

This implies that $\tilde{\rho} \leq \rho(1+\epsilon) < 1$ for $N \geq N(\varepsilon)$, and hence that equation (6.1) follows. To prove equation (6.2) we notice that from equation (6.3), for any $n \geq 1$ and $N \geq N(\epsilon/4)$

$$\|\pi P^n - \pi \tilde{P}^n\| \leq \sum_{i=0}^{n-1} \|\pi(\tilde{P}^{n-i-1}(P - \tilde{P})(P^i - \pi)(\cdot))\| \leq C\epsilon \sum_{i=0}^{n-1} \rho^i \leq \frac{C\epsilon}{1-\rho},$$

and since $\|\pi - \check{\pi}^N\| = \lim_{n \to \infty} \|\pi P^n - \tilde{\pi}^N \tilde{P}^n\|$. We conclude the proof by taking $N(\epsilon, \rho) = N(\varepsilon) \vee N(\epsilon/4)$. $\square$

THE PSEUDO-MARGINAL APPROACH 21**7. Examples: Reversible jumps.** In this section we illustrate the potential of the pseudo-likelihood approach to Monte Carlo computations developed in this paper in the context of reversible jump MCMC [3] algorithms (RJMCMC hereafter), which are well known for their difficult implementation. We start with a toy example, for which the true marginals are known exactly, hence providing a simple ground truth, and illustrate the interest of the approach in a scenario, which in our opinion reflects the difficulties encountered in practice when implementing RJMCMC in more realistic and difficult scenarios. We then move on to a more substantial example related to variable selection for generalized linear models. We first show how an apparently reasonable RJMCMC applied to a seemingly simple nested models selection problem can easily fail to produce reliable results and demonstrate how our methodology can easily circumvent this problem and render the algorithm much more reliable.

7.1. *Toy example.* We consider here a toy transdimensional target distribution defined on $\{1\} \times \mathbb{R} \cup \{2\} \times \mathbb{R}^2$,

$$\pi(\theta, z) = \mathbb{I}(\theta = 1)\tfrac{1}{4}\mathcal{N}(z; 0, 1)$$
$$+ \mathbb{I}(\theta = 2)\tfrac{3}{4}\mathcal{N}\left(z = \begin{bmatrix} x \\ y \end{bmatrix}; \begin{bmatrix} 0 \\ 0 \end{bmatrix}, \Sigma = \begin{bmatrix} 1 & -0.9 \\ -0.9 & 1 \end{bmatrix}\right).$$

In a Bayesian setup this would correspond to an inference problem for which two models $\mathcal{M}_1$ and $\mathcal{M}_2$ are considered to explain the data, the models being indexed by $\theta$. Obviously here $\pi(\theta = 1) = 1/4$ and $\pi(\theta = 2) = 3/4$. However in order to illustrate our methodology we develop here a reversible jump algorithm [3] to sample from this distribution, and compare our results with the exact distribution. A simple marginal ideal chain can be defined through the transition

$$(7.1) \quad P(\theta, \vartheta) = 1 \wedge \frac{\pi(\vartheta)}{\pi(\theta)}\mathbb{I}(\vartheta \neq \theta) + \mathbb{I}(\vartheta = \theta)\left[1 - 1 \wedge \frac{\pi(\vartheta')}{\pi(\theta)}\mathbb{I}(\vartheta' \neq \theta)\right],$$

that is in other words when in model $\mathcal{M}_1$ (resp. $\mathcal{M}_2$) we propose a jump to model $\mathcal{M}_2$ (resp. $\mathcal{M}_1$) with probability 1. This chain is obviously uniformly ergodic, which is often the case for finite discrete chains in practice. Now assume that we are given the algebraic expression for $\pi(\theta, z)$ up to a factor of $1/4$ and that, possibly with age, we fail to recognize 3 times a bivariate normal distribution for model $\mathcal{M}_2$. For simplicity, we will assume that we successfully recognize a univariate normal distribution for model $\mathcal{M}_1$. A standard approach in such situations consists of resorting to a RJMCMC algorithm that uses the (available to a constant) density $\pi(\theta, z)$ for model $\mathcal{M}_2$ and requires one to propose a bidimensional vector $z \sim \mathcal{Q}(\cdot)$ when attempting a move from model $\mathcal{M}_1$ to $\mathcal{M}_2$. Naturally the effectiveness of the algorithm will, as we shall see, highly depend on this proposal distribution whose choice might be far from obvious in more complex scenarios. In order to improve this basic RJMCMC algorithm, we investigate here a very simple strategy which relies on the pseudo-likelihood framework described earlier more sophisticated and efficient approaches are possible and currently being explored in other work. For any $\eta > 0$, let $\mathcal{N}_\eta(z; \mu, \Sigma_d)$ denote the truncated normal distribution such that $\mathcal{N}_\eta(z; \mu, \Sigma_d) \propto \mathcal{N}(z; \mu, \Sigma_d)$ whenever $(z - \mu)^{\mathrm{T}} \Sigma^{-1} (z - \mu) \leq \eta^2$ and $\mathcal{N}_\eta(z; \mu, \Sigma_d) = 0$ otherwise. For $\theta = 2$ we define $Q_\theta^N(Z) = \mathcal{Q}(z(1)) \prod_{i=2}^N \mathcal{N}_\eta(z(i); z(i-1) + \frac{\sigma^2}{2} \nabla_z \log \pi_\theta(z(i-1)), \sigma^2 I_2)$ for some $\sigma^2 \ll 1$, that is we use a form of discretization of the Langevin diffusion with drift $\frac{1}{2} \nabla_z \log \pi_\theta(z)$, and run a Markov chain with transition density $\Pi(z, \mathfrak{z}) = \mathcal{N}_\eta(\mathfrak{z}; z + \frac{\sigma^2}{2} \nabla_z \log \pi_\theta(z), \sigma^2)$ whose equilibrium distribution is an approximation of $\pi_\theta(z)$. The algorithm proceeds as the marginal algorithm described earlier in equation (7.1), except that for $\theta = 2$, $\pi(\theta)$ is replaced with the estimator

$$\tilde{\pi}^N(\theta) = \frac{1}{N} \left[ \frac{\pi(\theta, z(1))}{\mathcal{Q}(z(1))} \right.$$
(7.2)
$$\left. + \sum_{i=2}^N \frac{\pi(\theta, z(i))}{\mathcal{N}_\eta(z(i); z(i-1) + \sigma^2/2 \nabla_z \log \pi_\theta(z(i-1)), \sigma^2)} \right],$$

where $z(1), z(2), \ldots, z(N)$ are sampled according to $Q_\theta^N$ above. Note that the case $N = 1$ corresponds to the "standard" RJMCMC algorithm described above. For the purpose of illustration we took $\mathcal{Q}(z) = \mathcal{N}(z; [3, 3]^{\mathrm{T}}, I_2)$, which while being an obviously bad choice ensures irreducibility, and hence (in theory) convergence for the standard RJMCMC algorithm. We ran our algorithm for $N = 1, 5$ and $10$ for 450,000, 90,000 and 45,000 iterations, respectively, resulting in comparable computational efforts. The respective empirical expected acceptance probabilities were 0.0121, 0.5206 and 0.5056 (note that the theoretical expected acceptance probability in the stationary regime for the marginal algorithm is $1/4 + 3/4 \times 1/3 = 1/2$). In Figure 1 we present the "instantaneous" model probability estimators $1/k \sum_{i=1}^k \mathbb{I}(\theta_i = 1)$ and $1/k \sum_{i=1}^k \mathbb{I}(\theta_i = 2)$ as a function of the iterations $k$. Note the deceptive behavior observed for $N = 1$ which suggests that convergence has occurred.

7.2. *Application to variable selection in GLMs.* In this section we present an application of the pseudo-marginal principle to model selection in generalized linear models, and focus here more particularly on the logit link. More precisely we assume that we observe $M \geq 1$ realizations $(y_i, x_i)$ for $i \in \{1, \ldots, M\}$ of a random variable pair $(Y, X)$ taking values in $\{0, 1\} \times \mathbb{R}^k$ for some $k \geq 1$ and that the dependence between $Y$ and $X$ is characterized

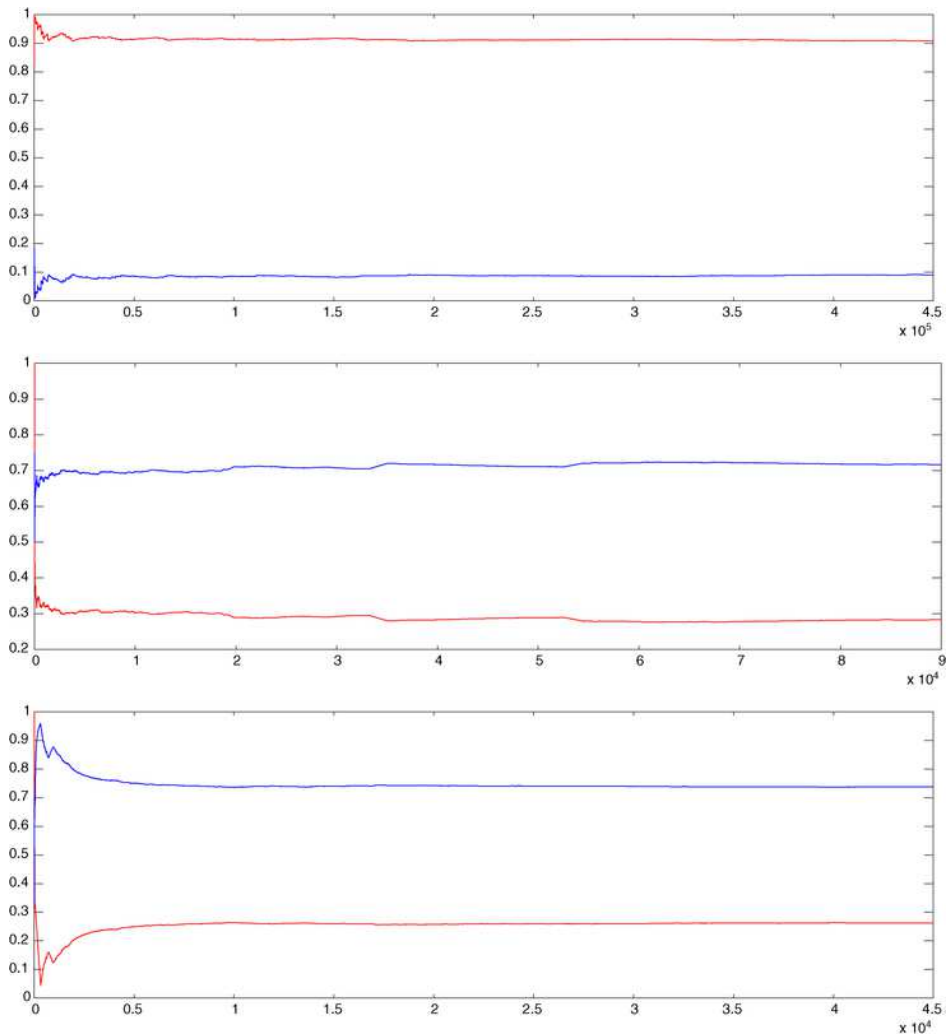

Fig. 1. *Instantaneous estimation of the model probabilities as a function of the iterations for $N = 1, 5, 10$ from top to bottom.*

by the conditional distribution $\mathbb{P}(Y = 1|X, z)$, assumed to satisfy

$$\log\left(\frac{\mathbb{P}(Y = 1|X, z)}{1 - \mathbb{P}(Y = 1|X, z)}\right) = Xz,$$

for a column vector $z \in \mathbb{R}^k$. Not all components $X(l)$, $l = 1, \ldots, k$, of $X$ (the "covariates") might be relevant to sparsely explain the data and as a result we might seek to compare models for which some of the components of $z$ are set to zero this is what we refer to as the model selection problem. In order to carry out inference in a Bayesian setup it is convenient to introduce indicator variables $\theta(l) \in \{0,1\}$ for $l = 1, \ldots, k$ such that covariate $X(l)$ is excluded whenever $\theta(l) = 0$. This allows us to index the $2^k - 1$ models (we exclude the model with no dependence) with the vector $\theta := (\theta(1), \theta(2), \ldots, \theta(k))$. Let $C$ be the $M \times k$ matrix whose $i$th row is $x_i$, we then denote $C_\theta$ the submatrix of $C$ that contains the columns $C$ for which $\theta(l) = 1$ and likewise for $z_\theta$ the subvector of $z$. It will be convenient to denote $\bar{\theta} = \sum_{l=1}^{k} \theta(l)$ the number of active covariates in model $\theta$. We ascribe prior distributions to $k$ and $z_\theta$ : $\Pr(\bar{\theta} = k) \propto \lambda^k/k!$ for some fixed $\lambda > 0$ and following [6] we set $z_\theta \sim \mathcal{N}(0, [C_\theta^T C_\theta/4M]^{-1})$ a priori. Denoting $p_i(z_\theta) := \mathbb{P}(y_i = 1|x_i, z_\theta)$ the joint posterior distribution is

$$\pi(\theta, z_\theta) \propto \prod_{i=1}^{N} p_i(z_\theta)^{\theta(i)}(1 - p_i(z_\theta))^{1-\theta(i)} \mathcal{N}(z_\theta; 0, [C_\theta^T C_\theta/4M]^{-1}) \lambda^{\bar{\theta}}/\bar{\theta}!.$$

Variable selection in a Bayesian context typically relies on the marginal posterior model probabilities $\pi(\theta)$, which are in the present situation intractable. One can for example resort to MCMC and this is the route followed up here.

The basis of our algorithm is a reversible jump MCMC algorithm for the marginal model which consists of a birth/death update. We first describe a *marginal algorithm* with transition $P(\theta, d\vartheta)$ which of course cannot be implemented. The pseudo-marginal algorithm will be a simple variation of this algorithm. Given a model $\theta$, $P(\theta, d\vartheta)$ can be described algorithmically as follows. With probability $1/2$, either:

- set $\theta^+ = \theta$ and if $\bar{\theta} < k$,
  1. choose uniformly among the $k - \bar{\theta}$ *nonactive* components' indexes, say $j$, and set $\theta^+(j) = 1$,
  2. set $\vartheta = \theta^+$ with probability

$$(7.3) \qquad 1 \wedge \frac{\pi(\theta^+)}{\pi(\theta)} \frac{1/(\bar{\theta}+1)}{1/(k-\bar{\theta})},$$

  otherwise $\vartheta = \theta$,

or

- set $\theta^- = \theta$ and if $\bar{\theta} > 0$,
  1. choose uniformly among the $\bar{\theta}$ *active* components' indexes, say $j$, and set $\theta^-(j) = 0$,
  2. set $\vartheta = \theta^-$ with probability

$$(7.4) \qquad 1 \wedge \frac{\pi(\theta^-)}{\pi(\theta)} \frac{1/(k-\bar{\theta}+1)}{1/\bar{\theta}},$$

  otherwise $\vartheta = \theta$.

This algorithm cannot be implemented, but it is nevertheless possible to implement a reversible jump algorithm on the joint distribution $\pi(\theta, z_\theta)$. A solution suggested in [2, 6], which we will refer to as the *standard RJ algorithm* here, can be understood as being a simple variation on the algorithm above, where in the birth move the additional sampling of a new coefficient from a distribution $\mathcal{Q}$ is required, resulting in the acceptance ratio

$$1 \wedge \frac{\pi(\theta^+, z_{\theta^+}^+)}{\pi(\theta, z_\theta)} \frac{1}{\mathcal{Q}(z^+(j))} \frac{1/(\bar{\theta}+1)}{1/(k-\bar{\theta})}.$$

In [6] it is suggested to use as a proposal distribution for $z^+(j)$ the marginal of the normal distribution with mean the maximum likelihood estimator $z_{1111}^{ML}$ of the saturated model and covariance the corresponding Hessian. Our pseudo-likelihood algorithm is very similar to the algorithm developed for the toy example in the previous section, that is it relies on a discretized Langevin diffusion, and consists formally of simply replacing $\pi(\theta^+)$ and $\pi(\theta^-)$ in the pseudo-code above with an estimator of the form equation (7.2) for $\theta \in \{\theta^+, \theta^-\}$. The following setup was considered. We generated artificial data from the logit model as follows. We chose $M = 50$, $k = 4$, a set of coefficients $z^* = [1\ 0.5\ -2\ 0.01]^\mathrm{T}$ and generated covariates as follows: with $Z_i \sim \mathcal{N}(0, I_M)$ we set $C_i = Z_i$ for $i = 1, 3, 4$ and $C_2 = 0.9 \times Z_1 + 0.1 \times Z_2$. This resulted into two correlated covariates, number 1 and 2. The maximum likelihood estimate for $z^*$ was found to be $z_{1111}^{ML} = [5.22445\ -3.71672\ -2.4011\ -0.587472]^\mathrm{T}$, suggesting (a) the presence of a main mode for the saturated model around this value, significantly different from the truth and (b) a mismatch between the modes and marginal modes of $\pi(1111, z_{1111})$ with those of $\pi(1011, z_{1011})$ ($z_{1011}^{ML} = [1.73253\ -2.30933\ -0.648927]^\mathrm{T}$) and $\pi(0111, z_{0111})$ ($z_{0111}^{ML} = [1.5968\ -1.98855\ -0.0922961]^\mathrm{T}$) which might result in poor mixing of the standard "birth–death" RJMCMC algorithm described above, a behavior likely to be reinforced here by the choice of the proposal distribution. This is confirmed by our simulation. In Figure 2 we present the estimated model probabilities (indexed by the decimal representation of $\theta$) for $N = 1, 5, 50, 100, 200$. Note that for the case $N = 1$ the birth/death move was complemented by 10 iterations of a within model one variable at a time random walk MH for each sweep. We observe the large discrepancy between the results for $N = 1, 5$ and the results for $N = 50, 100, 200$ the latter being in agreement. Note that this is despite the apparent convergence of estimators of the posterior inclusion probabilities $\mathbb{P}(\theta(j) = 1)$ for $j = 1, \ldots, k$ for the case $N = 1$ (Figure 3) and that the results obtained after 20,000, 10,000 and 5,000 for $N = 50, 100$ and 200 are much more reliable than for $N = 1$ after 1,000,000 iterations. The respective expected acceptance probabilities are given in the table below. Note that using the estimated model probabilities obtained for $N = 200$ one finds an acceptance rate of 0.29592. Finally, in

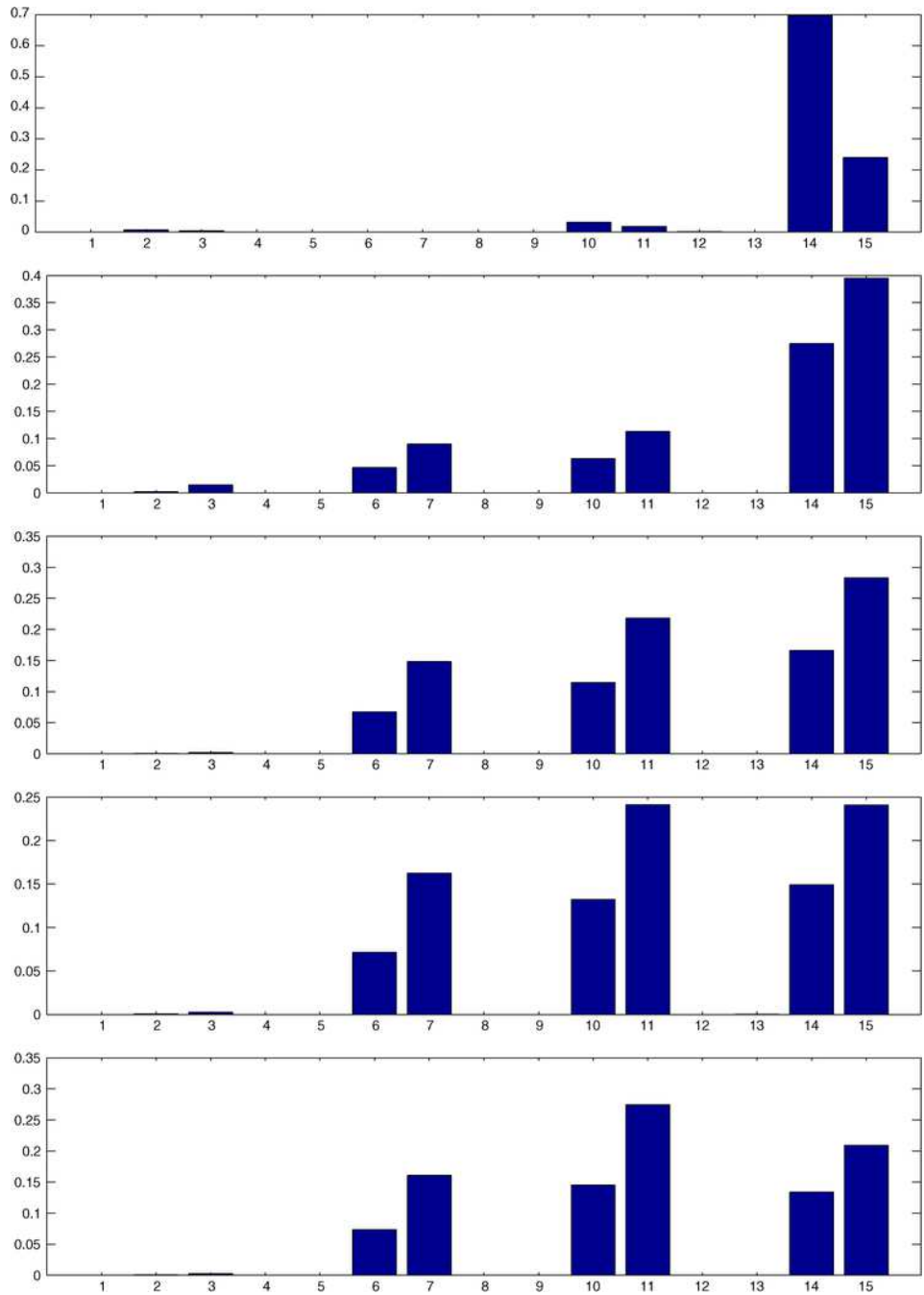

Fig. 2. *Model probabilities for $N = 1, 5, 50, 100, 200$ (from top to bottom).*

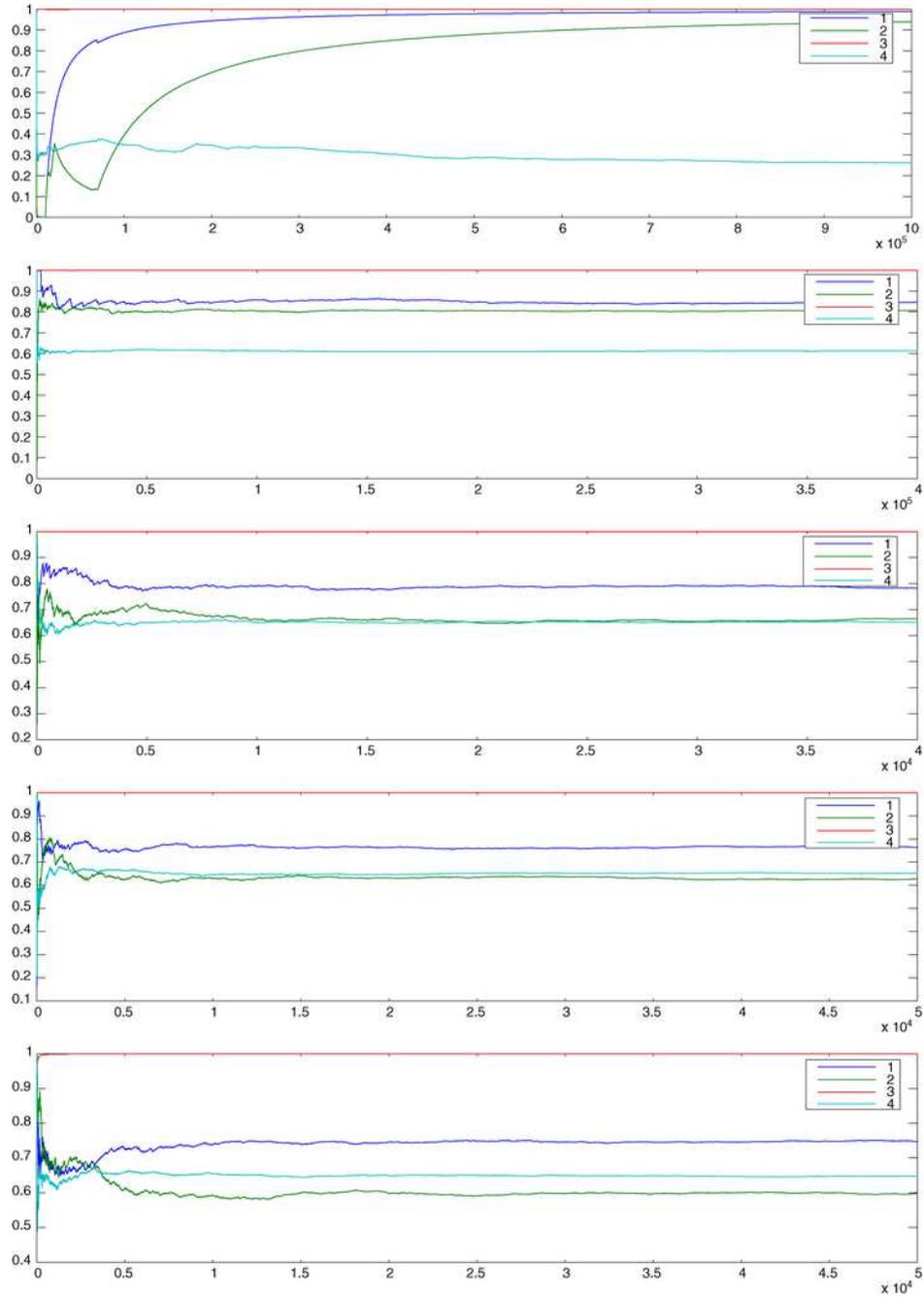

Fig. 3. *Instantaneous estimation of the model probabilities as a function of the iterations for $N = 1, 5, 50, 100, 200$.*

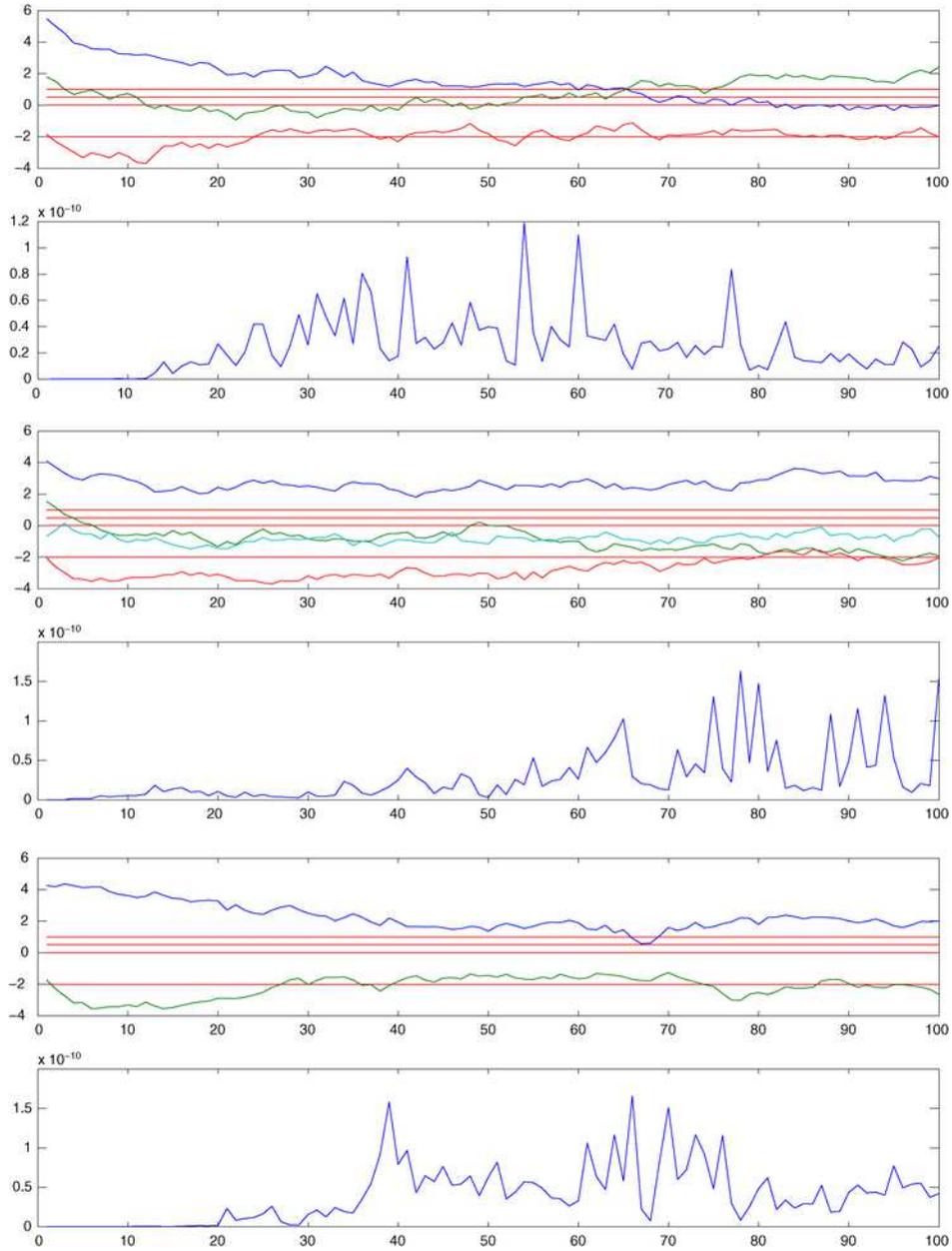

FIG. 4. *Snapshots of the transitions for $z_\theta$ between models for $N = 100$, together with their associated importance weights (the horizontal lines correspond to the true values of $z^*$).*

Figure 4 we present the trace of the $z$'s drawn with our discretized Langevin while attempting to jump between models, which together with the evolution of the values of the importance weights illustrates why our approach might be of interest.

| $N$/nb iter.$\times 1000$ | 1/1000 | 5/400 | 50/40 | 100/50 | 200/50 |
|---|---|---|---|---|---|
| accept prob. | 0.064293 | 0.16569 | 0.25885 | 0.28433 | 0.29371 |

**8. Conclusion.** The pseudo-marginal approach to stochastic simulation is a highly versatile methodology which has diverse potential applications in a variety of areas. The focus of this paper has been on some of the theoretical underpinnings of the method. Our main results describe ergodicity and uniform ergodicity of GIHM and its exact generalizations suggested in this paper, and we also give a comparison with an inexact variants, akin to MCWM. Empirical evidence in [1] and in the present paper in the context of reversible jumps for model selection in generalized linear models suggests that the methodology has considerable promise. Currently ongoing work confirms this.

## REFERENCES

[1] BEAUMONT, M. A. (2003). Estimation of population growth or decline in genetically monitored populations. *Genetics* **164** 1139–1160.
[2] DELLAPORTAS, P. and FORSTER, J. J. (1999). Markov chain Monte Carlo model determination for hierarchical and graphical log-linear models. *Biometrika* **86** 615–633. MR1723782
[3] GREEN, P. J. (1995). Reversible jump MCMC computation and Bayesian model determination. *Biometrika* **82** 711–732. MR1380810
[4] MEYN, S. P. and TWEEDIE, R. L. (1993). *Markov Chains and Stochastic Stability*. Springer, New York. MR1287609
[5] LIU, J. S., WONG, W. H. and KONG, A. (1994). Covariance structure of the Gibbs sampler with applications to the comparisons of estimators and augmentation schemes. *Biometrika* **81** 27–40. MR1279653
[6] NTZOUFRAS, I., DELLAPORTAS, P. and FORSTER, J. (2003). Bayesian variable and link determination for generalized linear models. *J. Statist. Plann. Inference* **111** 165–180. MR1955879
[7] O'NEILL, P. D., BALDING, D. J., BECKER, N. G., EEROLA, M. and MOLLISON, D. (2000). Analyzes of infectious disease data from houseing the expected value of ratios, hold outbreaks by Markov chain Monte Carlo methods. *Appl. Statist.* **49** 517–542. MR1824557
[8] ROBERT, C. P. and CASELLA, G. (2004). *Monte Carlo Statistical Methods*, 2nd ed. Springer, New York. MR2080278
[9] ROBERTS, G. O. and SAHU, S. K. (1997). Updating schemes, correlation structure, blocking and parameterisation for the Gibbs sampler. *J. Roy. Static. Soc. Ser. B* **59** 291–397. MR1440584
[10] ROBERTS, G. O. and TWEEDIE, R. (1996). Geometric convergence and central limit theorems for multidimensional Hastings and Metropolis algorithms. *Biometrika* **83** 95–110. MR1399158

[11] TIERNEY, L. (1998). A note on Metropolis–Hastings kernels for general state-spaces. *Ann. Appl. Probab.* **8** 1–9. MR1620401

SCHOOL OF MATHEMATICS  
UNIVERSITY OF BRISTOL  
BRISTOL, BS8 1TW  
UNITED KINGDOM  
E-MAIL: c.andrieu@bristol.ac.uk

DEPARTMENT OF STATISTICS  
UNIVERSITY OF WARWICK  
COVENTRY, CV4 7AL  
UNITED KINGDOM  
E-MAIL: Gareth.O.Roberts@warwick.ac.uk